\documentclass{article}[13pt]

\usepackage{fullpage}

\newtheorem{definition}{Definition} 
\newtheorem{theorem}{Theorem}
\newtheorem{proposition}{Proposition} 

\newtheorem{observation}{Observation}

\usepackage[all]{xy} \usepackage{amssymb}

\def\proof {Proof. \hspace{4pt}} \def\endofproof {\hfill $\Box$}

\def\gg {\mathfrak{g}} \def\hh {\mathfrak{h}}

\begin{document}

\title{AQFT from $n$-functorial QFT}

\author{Urs Schreiber\thanks{{\tt urs.schreiber@math.uni-hamburg.de}}}

\maketitle

\begin{abstract}
  There are essentially two different approaches to the
  axiomatization of quantum field theory (QFT): algebraic QFT, going back
  to Haag and Kastler, and
  functorial QFT, going back to Atiyah and Segal.
  More recently, based on ideas by Baez and Dolan,
   the latter is being refined to
  ``extended'' functorial QFT by Freed, Hopkins, Lurie and others.
  The first approach uses local nets of operator algebras
  which assign to each patch an algebra ``of observables'',
  the latter uses $n$-functors which assign to
  each patch a ``propagator of states''.

  In this note we present an observation about how these two
  axiom systems are naturally related: we demonstrate
  under mild assumptions that every 2-dimensional
  extended Minkowskian QFT 2-functor ("parallel surface transport")
  naturally yields a local net.
  This is obtained by postcomposing the propagation 2-functor
  with an operation that mimics the passage from the
  \emph{Schr{\"o}dinger picture} to the \emph{Heisenberg picture}
  in quantum mechanics.
  The argument has a straightforward generalization
  to general Lorentzian structure and higher dimensions.
\end{abstract}

\tableofcontents \newpage

\section{Introduction}

  Out of the numerous tools and concepts that physicists have used
  for the description of quantum field theory
  few are well defined beyond simple toy examples. Still, in many
  cases they ``work'', often with dramatic success. Axiomatizations
  of QFT attempt to extract from the ill-defined symbols that
  appear in the physics literature those properties which are actually
  being used in structural proofs.

  \begin{itemize}
    \item
     While the path integral
  itself usually is ill-defined, all that often matters is
  the assumption that it satisfies the gluing law
  \cite{Walker}. Taking
  this law as an axiom leads to the Atiyah-Segal formulation of
  functorial QFT.
   \item
    Similarly, while the products of physical field observables
    are usually ill-defined, all that often matters is the
    assumption that they satisfy the locality property
    \cite{BuchholzHaag}.
    Taking this as an axiom leads to the Haag-Kastler
    formulation of algebraic QFT.
  \end{itemize}

  The power of axiomatizations is that they lead to a more robust
  and clearer picture. The danger of axiomatizations is that they
  fail to capture important phenomena. Therefore it is
  especially important to understand how different axiomatizations
  of the same situation are related.

  \paragraph{AQFT: nets of local algebras.}
  Nets of local operator algebras have been introduced
  \cite{Haag} (see \cite{HalvorsonMueger} for a review)
  in order to formalize the concept of the algebra
  of \emph{local observables} in quantum field theory.
  One way to think of such a net is as a co-presheaf
  on a sub-category of open subsets
  of a given Lorentzian manifold $X$ with values in algebras.
  These co-presheaves
  are required to satisfy a couple of conditions
  (the first two mandatory, the third and fourth usually desired
  but sometimes dropped):
  \begin{enumerate}
    \item ({\bf isotony}) all co-restriction morphisms
      are required to be inclusions of sub-algebras -- this makes the
      co-presheaf a \emph{net};
    \item ({\bf locality/``microcausality''})
      the inclusions of two algebras assigned to
      two spacelike separated
      open subsets into the algebra assigned to a joint
      superset are required to commute with each other.
    \item
      ({\bf covariance})
       the net is covariant with respect to the action of a
       group $G$ on $X$ (for instance the Poincar{\'e}-group or the
       conformal group)
       if there is a family of algebra isomorphisms between
       the algebras assigned to any region and its image under the group
       action, compatible with the group product and the net structure.
    \item
      ({\bf time slice axiom}) the algebra of a subset is
      isomorphic to that assigned to any neighbourhood of any of its
      Cauchy surfaces.
  \end{enumerate}

  Out of the study of these structures a large subfield
  of mathematical physics has developed, which is
  equivalently addressed as \emph{algebraic quantum field theory},
  or as \emph{axiomatic quantum field theory} or as
  \emph{local quantum field theory}, but usually abbreviated
  as {\bf AQFT}. For a review of physical applications see
  \cite{whereweare}.

  \paragraph{FQFT: $n$-functorial cobordism representation.}
  Remarkably, all three of the terms -- algebraic, axiomatic, local
  -- would
  equally well describe what is probably the
  main alternative parallel development: the study of
  representations of cobordism categories, i.e.
  of functors from categories whose objects are
  $(d-1)$-dimensional manifolds and whose morphism
  are $d$-dimensional cobordisms between these to
  a category of vector spaces. An pedagogical introduction
  to this concept is in \cite{Baez}.

  Such functors have been introduced to formalize
  the concept of the quantum propagator acting
  on the space of quantum states and imagined to arise
  from an integral kernel given by a path integral.
  While this functorial approach did not receive a canonical
  name so far, here we shall refer to it as
  \emph{functorial quantum field theory} and
  abbreviate that as {\bf FQFT}.

  FQFT has most famously been studied in the context of
  \emph{topological} QFT, from which Atiyah originally deduced his
  sewing axioms \cite{Atiyah}. A review is \cite{Bartlett}.
  While topological FQFT is by far the most tractable and hence the
  best understood one,
  FQFT is not restricted to the topological case: equipping the
  cobordisms for instance with conformal structure yields
  conformal QFT, an observation which is the basis of Segal's
  functorial axiomatization of QFT \cite{Segal}.
  Restricting to 2-dimensional conformal cobordisms of genus 0
  this yields the axioms of vertex operator algebras
  \cite{Huang}, see \cite{Kong} for review and generalization.
  The result in \cite{FRS} can be regarded as providing examples
  for Segal's CFT axioms (though in that work Atiyah's formulation
  of the functoriality axiom is being referred to).

  Similarly, ordinary non-relativistic
  quantum mechanics ((1+0)-dimensional QFT) is about
  (monoidal) representations (i.e. functors to $\mathrm{Vect}$)
  of the (monoidal) category of 1-dimensional
  \emph{Riemannian} cobordisms \cite{ST}.
  Taking this point of view on ordinary quantum mechanics seriously
  leads to Abramsky-Coecke's \emph{categorical semantics of quantum
  protocols} \cite{AC}. See \cite{Coecke} for an overview.

  In this vein, here we shall be
  concerned with functors on cobordisms with
  \emph{pseudo}-Riemannian structures,
  and with flat Lorentzian structure (Minkowski structure)
  in particular.

  In \cite{FreedI,FreedII} it was suggested that the FQFT picture
  can and should be refined to an assignment of
  data of ``order $n$'' to codimension $n$ spaces
  for all $n$, such that this assignment respects all
  possible gluings. Formally this should mean that
  for $d$-dimensional quantum field theory the
  1-category of cobordisms is refined to
  a $d$-category of cobordisms
  \cite{ChengGurski,Verity} whose $k$-morphisms are $k$-dimensional
  cobordisms between $(k-1)$-dimensional cobordisms,
  and that one considers $d$-functors from this $d$-category
  to a suitable codomain $d$-category. Baez and Dolan
  began to draw the grand picture emerging here in \cite{BD},
  which was recently picked up by Hopkins and Lurie \cite{HL}.

  This extended $n$-functorial description of $d$-dimensional
  QFT is only beginning to be explored. First concrete
  descriptions of Chern-Simons and Wess-Zumino-Witten theory
  in this context appeared in \cite{FreedI,FreedII,ST}
  and in various talks given by Freed and Hopkins, aspects of 
  which have recently been made available as \cite{FreedIII}.
  Much progress has been made with understanding the
  extended FQFT of finite group Chern-Simons theory
  (Dijkgraaf-Witten theory) \cite{BW}.
  The general idea (for smooth $n$-groups) is currently best
  understood not for quantum but for ``classical''
  propagation, where it describes \emph{parallel transport}
  in $n$-bundles ($\simeq$ ($n-1$)-gerbes) with connection
  \cite{ndclecture,BS,SWI,SWII,SWIII}.

  But there are numerous indications that the picture is correct,
  useful and compelling. In \cite{FS} we shall demonstrate that
  the formulation of 2-dimensional CFT and 3-dimensional TFT
  appearing in \cite{FRS} (see \cite{FRSintro} for a review)
  is secretly a 2- and 3-FQFT of this form.

  \begin{table}[h]
  \begin{center}
  \begin{tabular}{c|c|c}
     {\bf names}
     &
     \begin{tabular}{c}
        algebraic QFT
        \\
        (also: axiomatic QFT,
        \\
        local QFT)
      \end{tabular}
      &
      functorial QFT
     \\
     \hline
     {\bf abbreviations}
     &
      {\bf AQFT}
      &
      {\bf FQFT}
      \\
      \hline
      &
      \multicolumn{2}{c}{assign}
      \\
      & algebras (of observables)
      & (time evolution) operators
      \\
      {\bf idea}
      & \multicolumn{2}{c}{to patches, compatible with}
      \\
      & inclusion & composition (gluing)
     \\
     \hline
     {\bf axioms due to}
     &
     Haag, Kastler
     &
     Atiyah, Segal
     \\
     \hline
     {\bf aspect of QFT}
     &
     \begin{tabular}{c}
        Heisenberg\\ picture
     \end{tabular}
     &
     \begin{tabular}{c}
        Schr{\"o}dinger \\ picture
     \end{tabular}
    \\
    \hline
    {\bf formal structure}
    &
    co-presheaf
    &
    \begin{tabular}{c}
      transport \\ $n$-functor
    \end{tabular}
     \\
     \hline
     {\bf
     \begin{tabular}{c}
        cartoon of
        \\
        domain structure
     \end{tabular}}
     &
    \raisebox{-20pt}{
    \xy
       (10,10)*{\bullet};
       (22,8)*{\bullet};
       (6,11)*{\bullet};
       (25,8)*{\bullet};
       \ar@{<-}_<{t} (0,20); (0,-5)
       \ar@{->}_>{x} (-5,0); (30,0)
       \ar@{-} (10,10); (15,15)
       \ar@{-} (15,15); (22,8)
       \ar@{-} (10,10); (17,3)
       \ar@{-} (17,3); (22,8)
       \ar@{-} (6,11); (14,19)
       \ar@{-} (14,19); (25,8)
       \ar@{-} (6,11); (17,0)
       \ar@{-} (17,0); (25,8)
    \endxy
    }
    &
    \raisebox{-20pt}{
    \xy
       (10,15)+(-7,-7)*{x}="s";
       (22,13)+(-7,-7)*{y}="t";
       (15,20)+(-7,-7)*{x'}="top";
       (17,8)+(-7,-7)*{}="bottom";
       (10,15)+(-2,-2)*{x'}="s2";
       (22,13)+(-2,-2)*{y'}="t2";
       (15,20)+(-2,-2)*{}="top2";
       (17,8)+(-2,-2)*{y}="bottom2";
       \ar@{<-}_<{t} (0,20); (0,-5)
       \ar@{->}_>{x} (-5,0); (30,0)
       \ar@{-} "s"; "top"
       \ar@{->} "top"; "t"
       \ar@{-} "s"; "bottom"
       \ar@{->} "bottom"; "t"
       \ar@{=>} "top"+(0,-3); "bottom"+(0,3)
       \ar@{-} "s2"; "top2"
       \ar@{->} "top2"; "t2"
       \ar@{-} "s2"; "bottom2"
       \ar@{->} "bottom2"; "t2"
       \ar@{=>} "top2"+(0,-3); "bottom2"+(0,3)
    \endxy
    }
    \\
    \hline
    {\bf relation}
    &
    \multicolumn{2}{c}{
     \begin{tabular}{c}
     \raisebox{18pt}{
      \xymatrix{
        \ar@{<-|}@/_2pc/[rrrr]|{\mbox{form endomorphism algebras}}
        &&&&
      }
      }
      \\
      \xymatrix@R=5pt{
    \makebox(35,38){
    \xy
       (10,15)+(0,-23)*{x}="s";
       (22,13)+(0,-23)*{y}="t";
       \ar@{-} "s"; (15,20)+(0,-23)
       \ar@{-} (15,20)+(0,-23); "t"
       \ar@{-} "s"; (17,8)+(0,-23)
       \ar@{-} (17,8)+(0,-23); "t"
    \endxy
     }
    \ar@{|->}[dd]^{A_Z}
     \\
      \hspace{30pt}
      &&&&
          Z
          \ar@{|->}[llll]
     \\
    \mathrm{End}
    \left(
    Z
    \left(
     \raisebox{0pt}{
    \makebox(30,15){
    \xy
       (10,15)+(0,-23)*{x}="s";
       (22,13)+(0,-23)*{y}="t";
       (17,8)+(0,-23)*{}="bottom";
       \ar@{-} "s"; "bottom"
       \ar@{->} "bottom"; "t"
    \endxy
    }
    }
    \right)
    \right)
      }
     \end{tabular}
    }
    \\
    \hline
    {\bf \begin{tabular}{c}
        main existing
        \\
        general theorems
     \end{tabular}}
     &
     \begin{tabular}{c}
       spin-statistics theorem,
       \\
       PCT theorem
     \end{tabular}
     &
     \begin{tabular}{l}
       results about
       \\
       topological invariants
     \end{tabular}
    \\
    \hline
    {\bf \begin{tabular}{c}
        main existing
        \\
        nontrivial examples
     \end{tabular}}
     &
     \begin{tabular}{c}
        chiral 2-d CFT
     \end{tabular}
     &
     \begin{tabular}{c}
       topological QFTs,
       \\
       full rational 2-d CFT
     \end{tabular}
  \end{tabular}
  \end{center}
   \caption{
     \label{summary table}
     {\bf The two approaches} to the axiomatization of
     quantum field theory together with their interpretation
     and relation as discussed here. The rectangular diagrams
     are explained in sections \ref{nets} and \ref{exFQFT}. The
     construction of the AQFT $\mathcal{A}_Z$ from the extended FQFT $Z$
     is our main point,
     described in section \ref{aqft from fqft}.}
  \end{table}

  \paragraph{The relation.}

  An obvious question, which does not seem to have
  been addressed before, is: \emph{What is the
   relation between the axioms of AQFT and FQFT?}

  Intuitively it is clear that the locality of local nets captures
  the same physical aspect as the $n$-functoriality of
  $n$-FQFTs does: that assignments to larger patches are
  already determined by the assignment to their pieces.
  But the nature of the assignments are different.
  We shall demonstrate that every FQFT determines an AQFT
  by postcomposing with the higher analog of the functor
  $$
    \mathrm{End} : \mathrm{Vect}_{\mathrm{iso}} \to
     \mathrm{Algebras}
  $$
  which sends each vector space to its algebra of endomorphisms
  and each isomorphism of vector space to the corresponding
  isomorphism of algebras.

  This functor is held in high esteem, if only implicitly so, in
  quantum mechanics, where it encodes the passage from what is
  called the \emph{Schr{\"o}dinger picture} to the
  \emph{Heisenberg picture} of quantum mechanics: given a
  unitary morphism
  of Hilbert spaces of the form
  $\xymatrix{E \ar[rr]^{e^{i t H}} && E}$ for $H$ some self-adjoint
 operator,
 which sends each element $\psi \in E$ to the element
 $e^{it H}\psi$,
 its image under the above functor is the isomorphism
 of endomorphism algebras
 $$
   \mathrm{End}
   :
   (\xymatrix{E \ar[rr]^{e^{i t H}} && E})
   \;
  \mapsto
   \;
   (\xymatrix{\mathrm{End}(E) \ar[rr]^{e^{i t H}\circ (-)\circ e^{-i t H}} &&
    \mathrm{End}(E)})
 $$
 which sends any operator $A$ on $E$ to $e^{it H}A e^{- it H}$.

  The situation is summarized in table \ref{summary table}.

  \clearpage

 \paragraph{Plan.}

  We start in section \ref{1dQFT} by discussing everything
  for the very simple case of 1-dimensional QFT (quantum mechanics),
  which should help to set the scene.
  Then in section \ref{nets} we quickly review those essentials of
  AQFT and in section \ref{exFQFT} those of FQFT which we need later
  on. Here we restrict to $d=2$ dimensions for ease of discussion.
  The generalization to higher dimensions is relatively obvious and
  straightforward, we briefly comment on that in section 
  \ref{higher dimensional QFT}.

  Our main definition is def. \ref{nets from functors}
  in section \ref{aqft from fqft}, which
  gives the prescription for turning an FQFT 2-functor into
  a 2-dimensional local net of algebras.
  Our main result is theorem \ref{consistency}, which states
  that this definition works.
  Theorem \ref{2funct construction} says that this construction extends to
  a 2-functor from the 2-category of FQFT 2-functors to the category
  of local nets,
  and, similarly, theorem
  \ref{covariance from equivariance}
  in section \ref{equivariance} says that
  the obvious notion
  of equivariance on FQFT induces the right notion of covariance in AQFT.

  We close by discussing some examples in section \ref{examples}
  and some further issues in section \ref{further issues}.

  \paragraph{2-categories.}
  See \cite{Leinster} for the basics of 2-categories
  and 2-functors between them. For the time being we can
  and will entirely restrict attention to \emph{strict} 2-categories and
  strict 2-functors between them. A review of all the basics of strict
  2-categories that we need here can be found for instance in the
  appendix of \cite{SWIII}. After we have established our construction
  for strict 2-categories the generalization to arbitrary weak
  2-categories is immediate.

  \paragraph{Acknowledgement.} 
  I am grateful to
  David Corfield, Christoph Schweigert, Zoran {\v S}koda, Jim Stasheff,
  Jamie Vicary
  and Konrad Waldorf for comments on earlier versions of this text, to
  Bruce Bartlett for discussion of aspects of some of the examples, to
  Maarten Bergvelt for discussion of relations with chiral nets and
  vertex operator algebras, to Jacques Distler for general
  discussion about AQFT and QFT, to Liang Kong for describing
  to me his work with Yi-Zhi Huang
  and to Peter Teichner for discussion of aspects at the beginning of section
  \ref{1dQFT}. 
  
  I had very useful discussion with Roberto Conti
  at an Oberwolfach CFT workshop in 2007, when I started thinking
  about the ideas presented here. Finally I heartily thank Pasquale Zito for
  a pleasant visit, for
  very useful discussion about his thesis and about Hopf spin chain models
  and for teaching me about asymptotic inclusion and pointing me to the
  relevant references.  

  This work was being
  completed while the author enjoyed a research fellowship at the
  \emph{Hausdorff Center for Mathematics} in Bonn.

\section{The situation for 1-dimensional QFT} \label{1dQFT}

  To put the following construction into perspective,
  it is useful to indicate what the transition from
  FQFT to AQFT that we are after looks like for the
  simple case where we are dealing with 1-dimensional
  quantum field theory, also known as quantum mechanics.

  \paragraph{Functorial quantum mechanics -- Schr{\"o}dinger picture.}

  There are some slight variations on the theme of how to think of
  ordinary quantum mechanics -- and in particular
  of possibly \emph{time dependent} quantum mechanics -- as a transport
  functor.
  These slight variations
  will have analogs also in higher dimensions, and hence are
  worth considering.

  Let $X = \mathbb{R}$ be the real line, thought of as the
  \emph{worldline} of a particle and in particular thought of
  as equipped with the obvious trivial Minkowski structure, which
  regards each vector as timelike. Let
  $P_1(X)$ be the category of homotopy classes of future-directed
  paths in $X$.
  Hence the objects of $P_1(\mathbb{R})$ are the points of $\mathbb{R}$
  and there is a unique morphism from $x$ to $y$ whenever
  $x \leq y$. In other words, $P_1(X)$ happens to be
  nothing but $\mathbb{R}$
  regarded as a poset.

  There is the closely related category,
  $1\mathrm{Cob}_{\mathrm{Riem}}$, whose objects are
  disjoint unions of points and whose morphisms are abstract 1-dimensional
  cobordisms
  equipped with a Riemannian structure. If we forget the monoidal
  structure on $1\mathrm{Cob}_{\mathrm{Riem}}$ (which is important,
  but not for our purposes here) and restrict it to just a single point,
  then we find
 $$
   1\mathrm{Cob}_{\mathrm{Riem}} \simeq \mathbf{B}\mathbb{R}_{0,+}
  =
  \left\{
    \xymatrix{
      \bullet \ar[r]^{t} & \bullet
    }
    \;|\;
    t \in [0,\infty)
  \right\}
   \,,
 $$
 where on the right we have the one-object category whose space of
 morphisms is the non-negative real half-line with composition given
 by addition of real numbers.
 There is a canonical projection functor
 $$
   \xymatrix{
     P_1(\mathbb{R})
     \ar@{->>}[r]
     &
     1\mathrm{Cob}_{\mathrm{Riem}}
   }
 $$
 which sends the path $\xymatrix{x \ar[r] & y}$ to the Riemannian
 cobordism $\xymatrix{\bullet \ar[rr]^{t=(y-x)} && \bullet}$ of the same
 lenght.

 Now, ordinary time-independent quantum mechanics is a functor
 $$
   Z : 1\mathrm{Cob}_{\mathrm{Riem}}
    \to \mathrm{Vect}_{\mathrm{isos}}
 $$
 which sends the single object of $1\mathrm{Cob}_{\mathrm{Riem}}$ to
 the \emph{space of states}, $E$, and sends the Riemannian cobordism of
 length $t$ to an automorphism
 $$
   Z
   :
   (\xymatrix{\bullet \ar[r]^t & \bullet})
   \mapsto
   (\xymatrix{E \ar[rr]^{\exp(i t H)} && E})
   \,,
 $$
 for $H$ some endomorphism of the complex vector space $E$ --
 the \emph{Hamiltonian}.
  Here we take $\mathrm{Vect}_{\mathrm{isos}}$ to be the category whose
  objects are vector space and whose endomorphisms are
  linear \emph{iso}morphisms.

 By the above, we can understand this as a functor
 on paths on the worldline,
 $P_1(\mathbb{R})$,
 which happens to factor through $\mathbf{B}\mathbb{R}_{0,+}$:
 $$
   \raisebox{40pt}{
   \xymatrix{
      P_1(\mathbb{R})
      \ar[r]
      \ar@{->>}[d]
      &
      \mathrm{Vect}_{\mathrm{isos}}
      \\
      \mathbf{B}\mathbb{R}_{0,+}
      \ar@{-}[r]^\simeq
      &
      1\mathrm{Cob}_{\mathrm{Riem}}
      \ar[u]_{Z}
   }
   }
   \,.
 $$
 Using the interpretation of such functors as vector bundles with
 connection \cite{SWI}, we can think of this as a vector bundle on
 the real line obtained from an $\mathbb{R}_{0,+}$-equivariant
 vector bundle over the point.

 A more general situation is obtained when one considers
 \emph{time dependent} quantum mechanics. Here the
 space of states and the Hamiltonian is allowed to change.
 There is then a 1-parameter family $t \mapsto E_t$ of
 spaces of states and $H$ is no longer necessarily constant.
This, then, is the case of a general functor $P_1(\mathbb{R}^2) \to
\mathrm{Vect}_{\mathrm{isos}}$: $$
  (\xymatrix{
     x \ar[r] & y
   })
   \mapsto
   (
    \xymatrix{
     E_x
     \ar[rrr]^{P \exp(i \int_x^y H(t)\, dt)}
     &&&
     E_y
    }
   )
   \,,
$$ where the expression on the right denotes the path-ordered exponential, which
is nothing but the parallel transport with respect to the connection 1-form $A =
H \, dt$. (More on that in section \ref{examples}.)

A slightly different but very similar concept plays an important role in
\cite{ST}, where quantum field theories \emph{over} a space $X$ are considered,
as functors from a category of cobordisms that come equipped with maps to $X$:
The category $1\mathrm{Cob}_{\mathrm{Riem}}(\mathbb{R})$ of cobordisms equipped
with a (smooth, say) map to the real line is not quite the same as
$P_1(\mathbb{R})$, but very similar. There is an obvious canonical functor $$
 \xymatrix{
   P_1(\mathbb{R})
   \ar[r]
   &
   1\mathrm{Cob}_{\mathrm{Riem}}(\mathbb{R})
} $$ which sends a path $\gamma$ in $\mathbb{R}$ to the Riemannian cobordism of
the same length equipped with the obvious map to $\mathbb{R}$ which coincides
with $\gamma$.

This way, from every ``1-dimensional QFT over $\mathbb{R}$'' in the sense of
\cite{ST} $$
  F : \xymatrix{
      1\mathrm{Cob}_{\mathrm{Riem}}(\mathbb{R})
      \to
      \mathrm{Vect}_{\mathrm{isos}}
   }
$$ one obtains an instance of ordinary time-dependent quantum mechanics by
pulling back to $P_1(\mathbb{R})$: $$
  \xymatrix{
    P_1(\mathbb{R})
    \ar[dr]
    \ar[rr]^{Z}
    &&
    \mathrm{Vect}_{\mathrm{isos}}
    \\
    &
    1\mathrm{Cob}_{\mathrm{Riem}}(\mathbb{R})
    \ar[ur]_F
  }
  \,.
$$ (In \cite{ST} Euclidean QFT is considered such that the morphisms assigned by
$Z$ are not in general invertible. While this is of no real relevance for the
point of the above discussion, notice that later on, when we pass from FQFT to
AQFT, we make crucial use of the fact that we assume FQFTs to assign invertible
time propagators.)

  Depending on the precise details, the functor $Z$ is
  usually demanded to factor through vector spaces
  with suitable extra structure. Topological vector spaces
  and Hilbert spaces are common choices. For our current purposes all
  such extra structure does not add anything to the aspects that
  we are interested in here and will be ignored until we come
  to concrete examples in section \ref{examples}.

\paragraph{Algebraic quantum mechanics -- Heisenberg picture.}

  Given such a functor $Z$, we can form for each point $x \in X$
  the \emph{endomorphism algebra} of the vector space, by sending
  $$
     x \mapsto \mathrm{End}(Z(x))
     \,.
  $$
  In the case that there is extra structure on our vector
  spaces we would demand suitable endomorphisms. In the case
  of Hilbert spaces one usually demands all endomorphisms to
  be \emph{bounded} operators.

  The endomorphism algebras thus obtained is known often as the
  \emph{algebra of observables}. In the present case,
  we would be tempted to associate this algebra at time $x$
  with the entire future of $x$.

  So let $S(X)$ be the category whose objects are
  open sets $ O_x := \{x' \in X| x' > x\}$ and whose
  morphisms are inclusions $O_x \subset O_{y}$ of open
  subsets.
  Of course, due to the simplicity of the present setup,
  $S(X)$ is canonically isomorphic to the opposite of
   $P_1(X)$ itself, hence
  is itself just the opposite catgeory of $\mathbb{R}$ regarded
  as a poset.
  But for the discussions to follow it is useful to
  think of $S(X)$ as a category of open subsets of $X$.

  The crucial point now is that sending spaces of states
  to their algebras of endomorphisms sends the functor
  $$
     Z : P_1(X) \to \mathrm{Vect}_{\mathrm{iso}}
  $$
  to a functor $\mathcal{A}_Z$ defined by
  $$
     \xymatrix{
         S(X) \ar[dr]_{Z} \ar[rr]^{A_Z} && \mathrm{Algebras}
         \\
         &
         \mathrm{Vect}_{\mathrm{iso}}
         \ar[ur]_{\mathrm{End}}
     }
     \,.
  $$
  The functor $\mathcal{A}_Z$ sends open subsets in $S(X)$ to the
  algebras of endomorphisms of the spaces of states
  sitting over their boundary, and it sends inclusions
  of open subsets to the inclusion of the algebras which
  is induced from using conjugation with the
  propagator that is assigned to the path connecting the
  respective boundaries. More precisely:
  $$
    \mathcal{A}_Z : (O_{y} \subset O_x)
    \mapsto
     (
        \xymatrix{
          \mathrm{End}(Z(y))
          \ar@{^{(}->}[rrr]^{Z(x \to y)^{-1}\circ(-) \circ Z(x \to y)}
          &&&
          \mathrm{End}(Z(x))
        }
     )
    \,.
  $$
  Of course this means that all inclusions of algebras
  here are actually isomorphisms.
  But this is again just due to the
  simplicity of the one-dimensional example. In conclusion, since there is no
  content in the locality axiom in 1 dimension, this means that
  $\mathcal{A}_Z$ is indeed a net of local monoids.

  It is this simple situation which we want to generalize
  from 1- to 2-dimensional QFT.

\section{Nets of local monoids} \label{nets}

 We start by considering a simple version of the relevant
 axioms of nets of local algebras on Minkowski space.
 Compare with section 2.1 of \cite{HalvorsonMueger}.
 Various refinements and generalizations are possible but add no
 further insight into the main point we want to make here.
 In particular, we shall ignore all extra structure that might be
 present on the algebras that appear below (such as them being
 $C^*$- or von-Neumann algebras) and even be content with regarding
 them just as \emph{monoids} (i.e. forgetting their vector space structure).
Our main point, that the inclusion and the locality axioms of local nets follow
from taking endomorphisms on $n$-functors, is entirely independent of all such
details. An interesting question is which extra structure on the $n$-functor
will induce which extra structure on the local nets. While this shall not
be our main concern here, the examples in section \ref{examples} give some
indications.

  So let $X = \mathbb{R}^2$ thought of as equipped with the
  standard
  Minkowski metric on $\mathbb{R}^2$.

  By a causal subset of $X$ we shall mean as usual the interior of the
  intersection of the future of one point with the past of another.

  \begin{figure}[h]
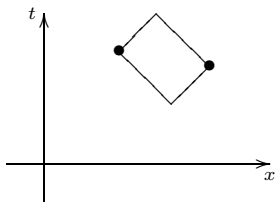

  $$
    \xy
       (10,15)*{\bullet};
       (22,13)*{\bullet};
       \ar@{<-}_<{t} (0,20); (0,-5)
       \ar@{->}_>{x} (-5,0); (30,0)
       \ar@{-} (10,15); (15,20)
       \ar@{-} (15,20); (22,13)
       \ar@{-} (10,15); (17,8)
       \ar@{-} (17,8); (22,13)
    \endxy
  $$
   \caption{
     A ``causal subset'' of 2-dimensional Minkwoski space
     is the interior of a rectangle all whose sides are
     lightlike. Such subsets are entirely fixed in particular
     by their left and right corners.
   }
  \end{figure}

  \begin{definition}
     \label{causal subsets}
     We denote by $S(X)$ the category whose objects are
     open causal subsets $V\subset X$ of $X$ and whose
     morphisms are inclusions $V \subset V'$.
  \end{definition}

  \begin{figure}[h]
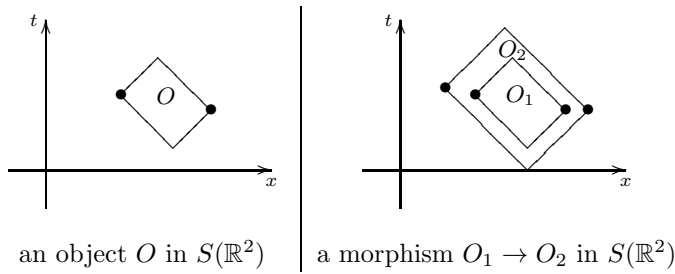

   \begin{center}
   \begin{tabular}{c|c}
   $$
    \xy
       (10,10)*{\bullet};
       (22,8)*{\bullet};
       (16,10)*{\mbox{\small $O$}};
       \ar@{<-}_<{t} (0,20); (0,-5)
       \ar@{->}_>{x} (-5,0); (30,0)
       \ar@{-} (10,10); (15,15)
       \ar@{-} (15,15); (22,8)
       \ar@{-} (10,10); (17,3)
       \ar@{-} (17,3); (22,8)
    \endxy
    $$
    &
    $$
    \xy
       (10,10)*{\bullet};
       (22,8)*{\bullet};
       (6,11)*{\bullet};
       (25,8)*{\bullet};
       (15,16)*{\mbox{\small $O_2$}};
       (16,10)*{\mbox{\small $O_1$}};
       \ar@{<-}_<{t} (0,20); (0,-5)
       \ar@{->}_>{x} (-5,0); (30,0)
       \ar@{-} (10,10); (15,15)
       \ar@{-} (15,15); (22,8)
       \ar@{-} (10,10); (17,3)
       \ar@{-} (17,3); (22,8)
       \ar@{-} (6,11); (14,19)
       \ar@{-} (14,19); (25,8)
       \ar@{-} (6,11); (17,0)
       \ar@{-} (17,0); (25,8)
    \endxy
    $$
    \\
    \\
    an object $O$ in $S(\mathbb{R}^2)$
    &
    a morphism $O_1 \to O_2$ in $S(\mathbb{R}^2)$
   \end{tabular}
   \end{center}
   \caption{
     The category $S(\mathbb{R}^2)$ of causal subsets of
     2-dimensional Minkowski space. Objects are causal
     subsets, morphisms are inclusions of these.
   }
  \end{figure}

  In order to concentrate just on the properties crucial for
  our argument, we
  shall now talk about nets of local \emph{monoids}
  (sets equipped with an associative and unital product).

  \begin{definition}
    Two objects $O_1$, $O_2$ in $S(X)$ are called
    spacelike separated if all pairs of points
    $(x_1,x_2) \in O_1 \times O_2$ are spacelike
    separated.
  \end{definition}

  \begin{figure}[h]
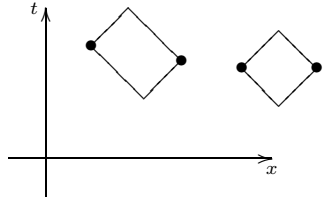


  $$
    \xy
       (10,15)+(-4,0)*{\bullet};
       (22,13)+(-4,0)*{\bullet};
       (26,12)*{\bullet};
       (36,12)*{\bullet};
       \ar@{<-}_<{t} (0,20); (0,-5)
       \ar@{->}_>{x} (-5,0); (30,0)
       \ar@{-} (10,15)+(-4,0); (15,20)+(-4,0)
       \ar@{-} (15,20)+(-4,0); (22,13)+(-4,0)
       \ar@{-} (10,15)+(-4,0); (17,8)+(-4,0)
       \ar@{-} (17,8)+(-4,0); (22,13)+(-4,0)
       \ar@{-} (26,12); (31,17)
       \ar@{-} (31,17); (36,12)
       \ar@{-} (26,12); (31,7)
       \ar@{-} (31,7); (36,12)
    \endxy
  $$
   \caption{
     Two spacelike separated causal subsets of $\mathbb{R}^2$.
   }
  \end{figure}

  \begin{definition}
    \label{netmonoids}
    A functor
    $$
      \mathcal{A} : S(\mathbb{R}^2) \to \mathrm{Monoids}
      \,,
    $$
    is a {\bf net} of monoids on 2-dimensional Minkwoski
    if it sends all morphisms in $S(\mathbb{R}^2)$ to injections
    (monomorphisms) of monoids.
    This is a {\bf net of local monoids}
    if for all spacelike separated $O_1,O_2 \subset O$
    the corresponding algebras commute with each other in $O$, i.e.
    $$
     [
       \mathcal{A}(O_1), \mathcal{A}(O_2)
     ]
       = 0
    $$
    as an identity in $\mathcal{A}(O)$.
      The net $\mathcal{A}$ is said to satisfy the {\bf time slice axiom}
    if for any region $O$, any Cauchy surface in $O$ and any collection
    of causal subset $\{O'_i \subset O\}$ covering the Cauchy surface
    we have
    $$
      \cup_i \mathcal{A}(O_i) = \mathcal{A}(O)
      \,,
    $$
    where the union is taken in $\mathcal{A}(O)$.
  \end{definition}
  Recall that a Cauchy surface of some region
  is a codimension 1 manifold such that
  all timelike or lightlike curves through that region cross it exactly once.
  In our case, Cauchy surfaces of a causal subset are all those curves
  through the subset which start at the left corner, monotonically
  move right, and end at the right corner.

 \begin{figure}[h]
  $$
    \xymatrix@R=10pt@C=10pt{
      &&&
      \ar@{~}[dddrrr]
      &&&&
      \\
      &&
      &&
      \ar@{-}[ddrr]
      \\
      &
      \ar@{-}[dr]
      &&
      \\
      x\ar@{-}[ur]
      \ar@{~}[uuurrr]
      \ar@{~}[dddrrr]
      \ar@{-}[ddrr]
      &&
      \ar@{-}[uurr]
      &&&& y
      \\
      &&&
      \ar@{-}[dr]
      \\
      &&\ar@{-}[ur]&&
      \ar@{-}[uurr]
      \\
      &&&
      \ar@{~}[uuurrr]
      &&
    }
  $$
  \caption{
    An inclusion $\{O'_i \subset O\}$ such that $\cup_i O'_i$
    contains Cauchy surfaces
    of $O$.
  }
 \end{figure}
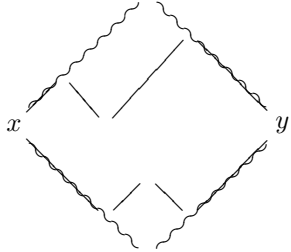

  Notice that a monoid (possibly an algebra) $A$ can be regarded as a one-object
  category
  $\mathbf{B}A := \left\lbrace \xymatrix{\bullet \ar[r]^a & \bullet} | a \in
  A\right\rbrace$
  (possibly enriched over vector spaces). As such, these monoids
  naturally form the
  2-category whose objects are monoids, whose morphisms are
  homomorphisms and whose 2-morphisms are intertwiners. See also
  appendix \ref{2veccanrep}.
  \begin{definition}
    We write
    $\mathbf{AQFT}(\mathbb{R}^2)$ for the
    sub-2-category
    of the 2-functor 2-category
    $2\mathrm{Funct}(S(\mathbb{R}^2), \mathrm{Cat})$
    whose objects are local nets $A$
    satisfying the time slice axiom, regarded as
    functors
    $$
      \xymatrix{
       S(\mathbb{R}^2)
       \ar[r]^{\mathcal{A}}
       &
       \mathrm{Monoids}
        \ar@{^{(}->}[r]^>>>>{\mathbf{B}(-)}
       &
       \mathrm{Cat}
      }
    $$
    taking values in one-object categories,
    whose morphisms are ordinary (as opposed to lax or pseudo)
    natural transformations between these, and whose 2-
    morphisms are modifications between those.
  \end{definition}

  \paragraph{Monoidal categories of endomorphisms of local nets.}
   From this it is immediate that
   for $\mathcal{A} \in \mathbf{AQFT}(\mathbb{R}^2)$
   the endomorphisms
   $
     \mathrm{End}_{\mathbf{AQFT}(\mathbb{R}^2)}(\mathcal{A})
   $
   form a monoidal category
    (since it arises from a one-object 2-category).
   This is the monoidal catgegory defined in
   definitions 8.1 and 8.5 in \cite{HalvorsonMueger} and proven there
   to be monoidal in proposition 8.30. The full subcategory
   $$
     \Delta(\mathcal{A}) \subset
     \mathrm{End}_{\mathbf{AQFT}(\mathbb{R}^2)}(\mathcal{A})
   $$
   of \emph{local}
   (meaning supported on some $O \in S(\mathbb{R}^2)$)
   and \emph{transportable} (meaning independent of support region up to
   isomorphism) endomorphisms is the main entity of interest in,
   and maybe in AQFT in general.
  The famous Doplicher-Roberts reconstruction theorem was motivated by
  the study of $\Delta(\mathcal{A})$. This is discussed in
  great detail in \cite{HalvorsonMueger}.

  \paragraph{Symmetries, covariance and equivariance.}
    Let $G$ be a group acting on $\mathbb{R}^2$ and preserving the
    causal set structure in that the action lifts to a functor
    $$
      g : S(\mathbb{R}^2) \to S(\mathbb{R}^2)
    $$
    for all $g \in G$. For $\mathcal{A}$ any local net we write
    $$
      g^*\mathcal{A} :
      \xymatrix{
        P_2(\mathbb{R}^2)
        \ar[r]^g
        &
        P_2(\mathbb{R}^2)
        \ar[r]^{\mathcal{A}}
        &
        \mathrm{Monoids}
      }
    $$
    for the pullback of the net along the action of $g\in G$.
    \begin{definition}
      \label{netequi}
      An equivariant structure on a local net $\mathcal{A}$
      is a choice of isomorphisms
      $$
        \xymatrix{
          \mathcal{A}
          \ar[r]^{r_g}
          &
          g^* \mathcal{A}
        }
      $$
      for all $g \in G$ such that for all $g_1,g_2 \in G$ we have
      $$
       \raisebox{40pt}{
        \xymatrix{
            & g_1^* \mathcal{A}
            \ar[dr]^{g_1^* r_{g_2}}
            \\
           \mathcal{A}
           \ar[ur]^{r_{g_1}}
           \ar[rr]^{g_1 g_2}
           &&
           (g_1 g_2)^* \mathcal{A}
        }
        }
        \,.
      $$
    \end{definition}

   \paragraph{Remark.} This is 1-categorical descent \cite{StreetDesc}
   along the nerve of the action groupoid $X//G$ of the category-valued
   presheaf $\mathrm{Funct}(S(-),\mathrm{Monoids})$.

   \paragraph{Remark.} In the AQFT literature this equivariant structure
   is often called a
   \emph{covariant} structure (for instance assumption 3 on p. 14
    of \cite{HalvorsonMueger}) and is often expressed in terms of
    the total algebra $\mathrm{colim}_{S(\mathbb{R}^2)} \mathcal{A}$
    (compare fact 5.10 on p. 41 of \cite{HalvorsonMueger}).

\section{Extended 2-dimensional Minkowskian FQFT} \label{exFQFT}

  Instead of regarding causal subsets as a category under
  inclusion of subsets, we can think of them as living in
  a 2-category under \emph{composition} (gluing).

  \begin{definition}
     \label{2paths in Minkowski}
     Let $P_2(\mathbb{R}^2)$ be the 2-category whose
     objects are the points of $\mathbb{R}^2$, whose
     morphisms are piecewise lightlike right-moving paths in $\mathbb{R}^2$
     and whose 2-morphisms are \emph{generated} from the
     closure of causal bigons
     $$
    
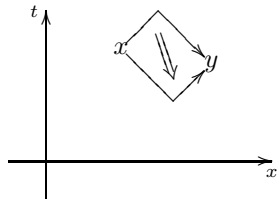
\xy
       (10,15)*{x}="s";
       (22,13)*{y}="t";
       (15,20)*{}="top";
       (17,8)*{}="bottom";
       \ar@{<-}_<{t} (0,20); (0,-5)
       \ar@{->}_>{x} (-5,0); (30,0)
       \ar@{-} "s"; "top"
       \ar@{->} "top"; "t"
       \ar@{-} "s"; "bottom"
       \ar@{->} "bottom"; "t"
       \ar@{=>} "top"+(0,-3); "bottom"+(0,3)
    \endxy
    $$
    regarded as 2-morphisms as indicated,
    under gluing along pieces of joint boundary.
    Composition is by gluing along pieces of joint boundary,
    in the obvious way.
  \end{definition}

  \begin{figure}[h]
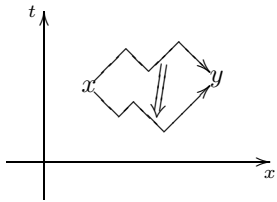

     $$
    \xy
       (6,10)*{x}="s";
       (23,11)*{y}="t";
       \ar@{<-}_<{t} (0,20); (0,-5)
       \ar@{->}_>{x} (-5,0); (30,0)
       \ar@{-} "s"; (11,15)
       \ar@{-} (11,15); (14,12)
       \ar@{-} (14,12); (18,16)
       \ar@{->} (18,16); "t"
       \ar@{-} "s"; (10,6)
       \ar@{-} (10,6); (12,8)
       \ar@{-} (12,8); (16,4)
       \ar@{->} (16,4); "t"
       \ar@{=>} (16,13); (15,6)
    \endxy
    $$
    \caption{
      A typical 2-morphism in $P_2(\mathbb{R}^2)$
    }
  \end{figure}

 \paragraph{Remark.} The restriction that 1-morphism have to go
 ``right'' and 2-morphisms ``downwards'' simplifies the discussion
 a bit but is otherwise of no real relevance. Various
 generalizations of $P_2(\mathbb{R}^2)$ can be considered without
 changing the substance of the following arguments.

  Just as with local nets, there are many variations of definitions
  of extended quantum field theories on 2-dimensional Minkowski
  space which one could consider. We choose to take the following
  simple definition. (Compare with the notion of parallel surface
  transport \cite{BS,SWII,SWIII}).

  \begin{definition}
    \label{eFQFT}
    For any 2-groupoid $C$,
  an extended FQFT on 2-dimensional Minkowski space is a 2-functor
  $$
    Z : P_2(\mathbb{R}^2) \to C
    \,.
  $$
  We write $\mathbf{FQFT}(\mathbb{R}^2,C) :=
  2\mathrm{Funct}(P_2(\mathbb{R}^2),C)$
  for the 2-functor 2-category and
  $\mathbf{FQFT}_{\mathrm{isos}}(\mathbb{R}^2,C)$
  for the maximal strict 2-groupoid inside it.
  \end{definition}

  In concrete application $C$ will usually be a 2-category
  of 2-vector spaces (which in general is not strict),
  as for instance those whose objects are (von Neumann) algebras,
  whose morphisms are bimodules over these, and whose 2-morphisms
  are bimodule homomorphisms \cite{ST}. We will see such an
  example in section \ref{examples} based on some constructions
  summarized in appendix \ref{2veccanrep}.

  But for the moment we do not need to make any concrete choice
  concerning $C$.
  The only necessary requirement for the following
  is actually that the 2-morphisms
  in $C$ all be invertible and that horizontal composition
  by the images of the 1-morphisms under $Z$ is injective.

\paragraph{Equivariant structures.} Let $G$ be a group acting by diffeomorphisms
on $\mathbb{R}^2$ which respects causal subsets in that the action extends to a
functor $$
  g : S(\mathbb{R}^2) \to S(\mathbb{R}^2)
$$ for all $g\in G$. There is a canonical notion of what it means for a
2-functor $Z : P_2(\mathbb{R}^2) \to C$ to be \emph{equivariant} with respect to
this action \cite{SWIII,equi,ndclecture}: for $g \in G$ denote by $$
  g^* Z :
  \xymatrix{
    P_2(\mathbb{R}^2)
    \ar[r]^{g}
    &
    P_2(\mathbb{R}^2)
    \ar[r]^Z
    &
    C
  }
$$ the pullback of $Z$ along the diffeomorphism $G$.

\begin{definition}[equivariance of 2-functors] \label{equivariance of
2-functors} A $G$-equivariant structure on $Z$ is choice of isomorphisms $f_g$
of 2-functors (i.e. strictly invertible pseudonatural transformations) $$
  \xymatrix{
    Z
    \ar[r]^{f_g}_{\simeq}
    &
    g^* Z
  }
$$ for all $g \in G$, and a choice for all $g_1,g_2 \in G$ of invertible
2-morphisms (i.e. modifications of pseudonatural transformations) $$
  \xymatrix{
    &
    g_1^* Z
    \ar[dr]^{g_1^* f_{g_2}}
    \\
    Z
    \ar[ur]^{f_{g_1}}
    \ar[rr]_{f_{g_1 g_2}}^{\ }="t"
    &&
    (g_1 g_2)^* Z
    \ar@{=>}^{F_{g_1,g_2}}_{\simeq} "t"+(0,8); "t"
  }
$$ such that for all $g_1,g_2,g_3 \in G$ the tetrahedra 2-commute: $$
  \raisebox{35pt}{
  \xymatrix{
    g_1^*Z
    \ar[rr]^{g_1^* f_{g_2}}
    &&
    (g_1 g_2)^* Z
    \ar[dd]^{(g_1 g_2)^* f_{g_3}}
    \\
    \\
    Z
    \ar[uu]^{f_{g_1}}
    \ar[rr]_{f_{g_1 g_2 g_3}}^{\ }="b"
    \ar[uurr]|{f_{g_1 g_2}}^{\ }="a"
    &&
    (g_1 g_2 g_3)^* Z
    \ar@{=>}^{F_{g_1,g_2}} "a"+(-8,8); "a"
    \ar@{=>}^{F_{g_1 g_2, g_3}} "b"+(4,8); "b"
  }
  }
  \hspace{10pt}
  =
  \hspace{10pt}
  \raisebox{35pt}{
  \xymatrix{
    g_1^*Z
    \ar[rr]^{g_1^* f_{g_2}}
    \ar[ddrr]|{g_1^* f_{g_2 g_3}}^{\ }="a"
    &&
    (g_1 g_2)^* Z
    \ar[dd]^{(g_1 g_2)^* f_{g_3}}
    \\
    \\
    Z
    \ar[uu]^{f_{g_1}}
    \ar[rr]_{f_{g_1 g_2 g_3}}^{\ }="b"
    &&
    (g_1 g_2 g_3)^* Z
    \ar@{=>}^{F_{g_2,g_3}} "a"+(8,8); "a"
    \ar@{=>}_{F_{g_1,g_2 g_3}} "b"+(-4,8); "b"
  }
  }
  \,.
$$ \end{definition}

\paragraph{Remark.} In the case that $G$ acts freely, this is nothing but
2-categorical descent \cite{StreetDesc} along $Y := (\xymatrix{X \ar@{->>}[r] &
X/G})$ with coefficients in the 2-category-valued presheaf
$2\mathrm{Funct}(P_2(-),C)$ \cite{ndclecture}. If $G$ does not act freely it is
descent with respect to the nerve of the action groupoid of $G$.

\section{The main point: AQFT from extended FQFT} \label{aqft from fqft}

  We define a map from FQFTs in the sense of definition
  \ref{eFQFT} to AQFTs in the sense of definition \ref{netmonoids}
  and demonstrate, theorem \ref{consistency}, that
  it indeed sends 2-functors to local nets of monoids satisfying the
  time slice axiom. Then we observe, theorem \ref{2funct construction},
  that this construction extends to
  a 2-functor from FQFTs to AQFTs on $\mathbb{R}^2$.

  \begin{definition}
    \label{nets from functors}
    Given any extended 2-dimensional FQFT, i.e. a 2-functor
    $$
      Z : P_2(\mathbb{R}^2) \to C
    $$
    we define a functor
    $$
      \mathcal{A}_Z : S(\mathbb{R}^2) \to \mathrm{Monoids}
      \,.
    $$
    On objects it it acts as
    $$
      \mathcal{A}_Z :
    \left(
    \raisebox{-20pt}{
    \xy
       (10,15)*{x}="s";
       (22,13)*{y}="t";
       \ar@{-} "s"; (15,20)
       \ar@{-} (15,20); "t"
       \ar@{-}|\gamma "s"; (17,8)
       \ar@{-} (17,8); "t"
    \endxy
     }
    \right)
    \hspace{5pt}
     \mapsto
    \hspace{5pt}
    \mathrm{End}_C
    \left(
    Z
    \left(
     \raisebox{-15pt}{
    \xy
       (10,15)*{x}="s";
       (22,13)*{y}="t";
       (17,8)*{}="bottom";
       \ar@{-}|\gamma "s"; "bottom"
       \ar@{->} "bottom"; "t"
    \endxy
    }
    \right)
    \right)
    \,,
    $$
    where on the right we form the monoid of 2-endomorphism $a$
    in $C$ on the 1-morphism $Z(x \stackrel{\gamma}{\to} y)$ in $C$
    that is the past boundary of $O_{x,y}$,
    $$
      \xymatrix{
        Z(x)
        \ar@/^2pc/[rr]^{Z(x\stackrel{\gamma}{\to} y)}_{\ }="s"
        \ar@/_2pc/[rr]_{Z(x\stackrel{\gamma}{\to} y)}^{\ }="t"
        &&
        Z(y)
        \ar@{=>}^a "s"; "t"
      }
      \,.
    $$
    On morphisms $\mathcal{A}_Z$ is defined to act as follows.
  \end{definition}
For any inclusion $O_{x',y'} \subset O_{x,y} \in S(\mathbb{R}^2)$ $$
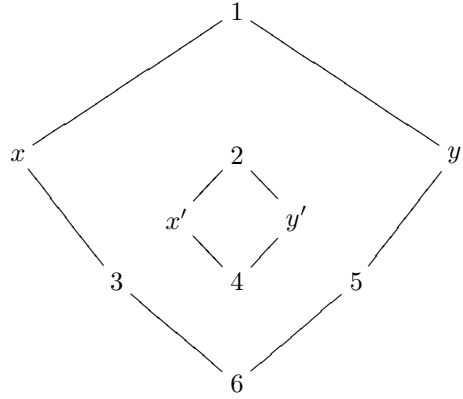
\xymatrix@C=10pt@R=10pt{
   &&&& 1
   \ar@{-}[dddrrrr]
   \\
   \\
   \\
   x
   \ar@{-}[ddrr]
   \ar@{-}[uuurrrr]
   &&&&
   2
   \ar@{-}[dr]
   &&&&
   y
   \\
   && & x' \ar@{-}[ur] \ar@{-}[dr]
    && y'
   %\ar@{-}[dr]
   &&
   \\
   && 3
   \ar@{-}[ddrr]
   %\ar@{-}[ur]
   && 4
    \ar@{-}[ur]
    %\ar@{=>}[dd]^f
   && 5
   \ar@{-}[uurr]
   \\
   \\
   &&&& 6
   \ar@{-}[uurr]
} $$ (the numbers here and in the following are just labels for various points
in order to help us navigate these diagrams) we form the pasting diagram $$
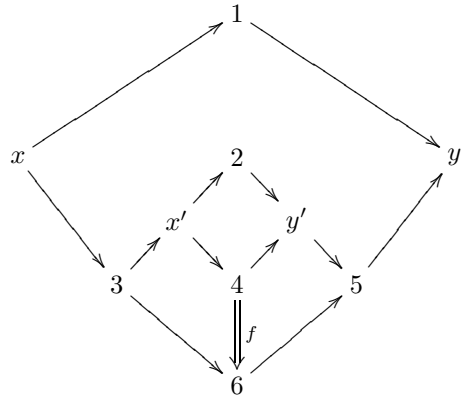
\xymatrix@C=10pt@R=10pt{
   &&&& 1
   \ar[dddrrrr]
   \\
   \\
   \\
   x
   \ar[ddrr]
   \ar[uuurrrr]
   &&&&
   2
   \ar[dr]
   &&&&
   y
   \\
   && & x' \ar[ur] \ar[dr]
    && y'
   \ar[dr]
   &&
   \\
   && 3
   \ar[ddrr]
   \ar[ur]
   && 4
    \ar[ur]
    \ar@{=>}[dd]^f
   && 5
   \ar[uurr]
   \\
   \\
   &&&& 6
   \ar[uurr]
} $$ in $P_2(\mathbb{R}^2)$. Here the obvious projections along light-like
directions (for instance from $x'$ onto $x \to 6$ yielding 3) is used. It is at
this point that the light-cone structure crucially enters the construction.

Let $f'$ be the 2-morphism obtained by whiskering (= horizontal composition with
identity 2-morphisms) the indicated 2-morphism $f$ with the 1-morphisms $x \to
3$ and $5 \to y$. $$
 f' := \hspace{4pt}
 \raisebox{34pt}
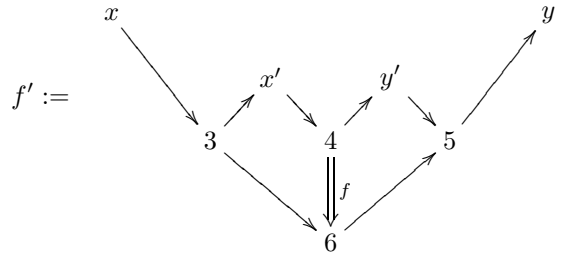
{
\xymatrix@C=10pt@R=10pt{
   x
   \ar[ddrr]
   &&&&
   &&&&
   y
   \\
   && & x'  \ar[dr]
    && y'
   \ar[dr]
   &&
   \\
   && 3
   \ar[ddrr]
   \ar[ur]
   && 4
    \ar[ur]
    \ar@{=>}[dd]^f
   && 5
   \ar[uurr]
   \\
   \\
   &&&& 6
   \ar[uurr]
}
 }
  \,.
$$ For any $a \in \mathrm{End}_C{Z(x',4,y')}$, 
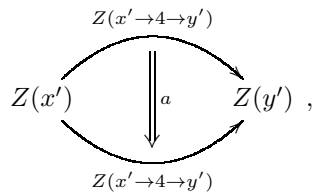
$$
  \xymatrix{
    Z(x')
    \ar@/^2pc/[rr]^{Z(x' \to 4 \to y')}_{\ }="s"
    \ar@/_2pc/[rr]_{Z(x' \to 4 \to y')}^{\ }="t"
    &&
    Z(y')
    \ar@{=>}^a "s"; "t"
  }
  \,,
$$ let $a'$ be the corresponding re-whiskering by $Z(x,3,x')$ from the left and
by $Z(y',5,y)$ from the right: $$
  \xymatrix{
    Z(x)
    \ar@/^2pc/[rr]^{Z(x\to 3 \to x' \to 4 \to y'  \to 5 \to y)}_{\ }="s"
    \ar@/_2pc/[rr]_{Z(x\to 3 \to x' \to 4 \to y'  \to 5 \to y)}^{\ }="t"
    &&
    Z(y)
    \ar@{=>}^{a'} "s"; "t"
  }
  \hspace{7pt}
    :=
  \hspace{7pt}
  \xymatrix{
    Z(x)
    \ar[rrr]^{Z(x \to 3 \to x')}
    &&&
    Z(x')
    \ar@/^2pc/[rr]^{Z(x' \to 4 \to y')}_{\ }="s"
    \ar@/_2pc/[rr]_{Z(x' \to 4 \to y')}^{\ }="t"
    &&
    Z(y')
    \ar[rrr]^{Z(y' \to 5 \to y)}
    &&&
    Z(y)
    \ar@{=>}^a "s"; "t"
  }
  \,,
$$ Then we obtain an injection $$
  \xymatrix{
     \mathrm{End}_C(Z(x',4,y'))
     \ar@{^{(}->}[rr]
     &&
     \mathrm{End}_C(Z(x,3,6,5,y))
   }
$$ by setting $$
  a \mapsto Z(f') \circ a' \circ Z(f')^{-1}
  \,,
$$ i.e. $$
  \xymatrix{
    Z(x')
    \ar@/^2pc/[rr]^{Z(x' \to 4 \to y')}_{\ }="s"
    \ar@/_2pc/[rr]_{Z(x' \to 4 \to y')}^{\ }="t"
    &&
    Z(y')
    \ar@{=>}^a "s"; "t"
  }
  \hspace{7pt}
    \mapsto
  \hspace{7pt}
  \xymatrix{
    Z(x)
    \ar[rrr]^{Z(x \to 3 \to x')}
    \ar@/^5pc/[rrrrrrrr]^{Z(x \to 3 \to 6 \to 5 \to y)}_{\ }="s2"
    \ar@/_5pc/[rrrrrrrr]_{Z(x \to 3 \to 6 \to 5 \to y)}^{\ }="t3"
    &&&
    Z(x')
    \ar@/^2pc/[rr]|{Z(x' \to 4 \to y')}_{\ }="s"^{\ }="t2"
    \ar@/_2pc/[rr]|{Z(x' \to 4 \to y')}^{\ }="t"_{\ }="s3"
    &&
    Z(y')
    \ar[rrr]^{Z(y' \to 5 \to y)}
    &&&
    Z(y)
    \ar@{=>}^a "s"; "t"
    \ar@{=>}^{Z(f')^{-1}} "s2"; "t2"
    \ar@{=>}^{Z(f')} "s3"; "t3"
  }
  \,.
$$

\paragraph{Remark.}
 Notice that this prescription is essentially nothing but the
 one we described already for the 1-dimensional case
 in section \ref{1dQFT}:
 to open subsets we assign the endomorphism algebra of the
 space of states assigned to one part of their boundary.
 To an inclusion of open subsets we then assign the
  inclusion of such algebras obtained by
  \emph{parallel transporting}
  the algebra of the inner set into the algebra of the
  outer set using conjugation with the propagators
  that the 2-functor assigns to 2-morphisms in $P_2(\mathbb{R}^2)$.
  The difference to the 1-dimensional case here is that this
  conjugation operation involves some (the obvious) re-whiskering.
  We will see that it is essentially this re-whiskering
  and the exchange law in 2-categories which
  lead to the locality of the net of monoids obtained this way.

  \begin{figure}[h]
    $$
      \xymatrix{
        a
        \ar@/^2pc/[rr]^{f_1}_{\ }="a1"
        \ar[rr]|{f_2}^{\ }="a2"_{\ }="a22"
        \ar@/_2pc/[rr]_{f_3}^{\ }="a3"
        &&
        b
        \ar@{}[r]|\cdot
        &
        b
        \ar@/^2pc/[rr]^{f'_1}_{\ }="b1"
        \ar[rr]|{f'_2}^{\ }="b2"_{\ }="b22"
        \ar@/_2pc/[rr]_{f'_3}^{\ }="b3"
        &&
        c
        \ar@{=>}^{F_1} "a1"; "a2"
        \ar@{=>}^{F_2} "a22"; "a3"
        \ar@{=>}^{F'_1} "b1"; "b2"
        \ar@{=>}^{F'_2} "b22"; "b3"
      }
      \hspace{5pt}
       =
      \hspace{5pt}
      \raisebox{14pt}{
      \xymatrix@R=0pt{
        a
        \ar@/^2pc/[rr]^{f_1}_{\ }="a1"
        \ar[rr]|{f_2}^{\ }="a2"
        &&
        b
        \ar@/^2pc/[rr]^{f'_1}_{\ }="b1"
        \ar[rr]|{f'_2}^{\ }="b2"
        &&
        c
        \\
        && \circ
        \\
        a
        \ar[rr]|{f_2}_{\ }="a22"
        \ar@/_2pc/[rr]_{f_3}^{\ }="a3"
        &&
        b
        \ar[rr]|{f'_2}_{\ }="b22"
        \ar@/_2pc/[rr]_{f'_3}^{\ }="b3"
        &&
        c
        \ar@{=>}^{F'_1} "b1"; "b2"
        \ar@{=>}^{F'_2} "b22"; "b3"
        \ar@{=>}^{F_1} "a1"; "a2"
        \ar@{=>}^{F'_1} "a22"; "a3"
      }
      }
    $$
    \caption{
     The exchange law in 2-categories, which is the functoriality of
     horizontal composition on the Hom-categories, says that the 2-dimensional
     order of composition of 2-morphisms is irrelevant.
    }
  \end{figure}

Now we come to our main point.

\begin{theorem}
  \label{consistency}
  The functor $\mathcal{A}_Z$ is a net of local monoids
  satisfying the time slice axiom.
\end{theorem} \proof
  We need to demonstrate three things
  \begin{enumerate}
     \item
       that the above assignment is functorial;
     \item
       that the above assignment satisfies the locality
       axiom;
     \item
       that the above assignment satisfies the time slice
       axiom.
  \end{enumerate}
  The third property is immediate from the construction.
  The first two properties turn out to be a direct consequence of
  2-functoriality of $Z$ and the exchange law in
  2-categories.

  To see functoriality, consider a chain of inclusions
  $$
    
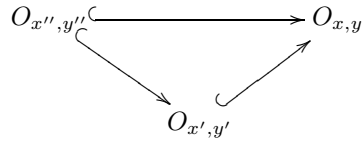
\xymatrix{
      O_{x'',y''}
      \ar@{^{(}->}[rr]
      \ar@{^{(}->}[dr]
      &&
      O_{x,y}
      \\
      &
      O_{x',y'}
      \ar@{^{(}->}[ur]
    }
  $$
  in $S(\mathbb{R}^2)$
  and the corresponding pasting diagram
  $$
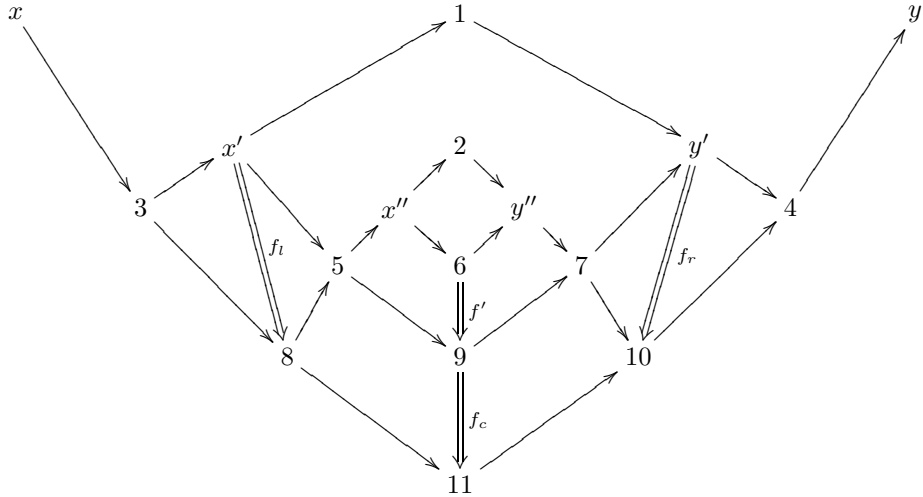
\xymatrix@C=8pt@R=8pt{
   x
   \ar[ddddrrr]
   &&
   &&&&&&& 1
   \ar[dddrrrr]
   &&&&&&& &&
   y
   \\
   \\
   \\
   &&&&& x'
   \ar[ddrr]
   \ar[uuurrrr]
   \ar@{=>}[ddddr]^{f_l}
   &&&&
   2
   \ar[dr]
   &&&&
   y'
   \ar[drr]
   \ar@{=>}[ddddl]^{f_r}
   &&&
   \\
   &&&   3
   \ar[urr]
   \ar[dddrrr]
   &&&& & x'' \ar[ur] \ar[dr]
    && y''
   \ar[dr]
   &&
   &&&
   4
   \ar[uuuurrr]
   \\
   &&&&
   &&& 5
   \ar[ddrr]
   \ar[ur]
   && 6
    \ar[ur]
    \ar@{=>}[dd]^{f'}
   && 7
   \ar[uurr]
   \ar[ddr]
   \\
   \\
   &&&&&&
      8
   \ar[dddrrr]
   \ar[uur]
   &&&
   9
   \ar[uurr]
   \ar@{=>}[ddd]^{f_c}
   &&&
   10
   \ar[uuurrr]
   \\
   \\
   \\
   &&&&&&&&& 11
   \ar[uuurrr]
} $$ in $P_2(\mathbb{R}^2)$. The composite inclusion $$
  \mathrm{End}_C(Z(x'' \to 6 \to y''))
  \hookrightarrow
  \mathrm{End}_C(Z(x' \to 5 \to 9 \to 7 \to y'))
  \hookrightarrow
  \mathrm{End}_C(Z(x \to 3 \to 8 \to 11 \to 10 \to 4 \to y))
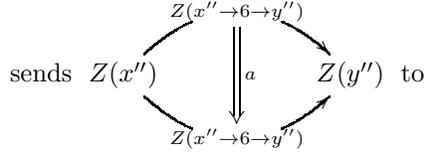
$$ sends $
  \xymatrix{
    Z(x'')
    \ar@/^2pc/[rr]|{Z(x'' \to 6 \to y'')}_{\ }="s"
    \ar@/_2pc/[rr]|{Z(x'' \to 6 \to y'')}^{\ }="t"
    &&
    Z(y'')
    \ar@{=>}^a "s"; "t"
  }
$ to 
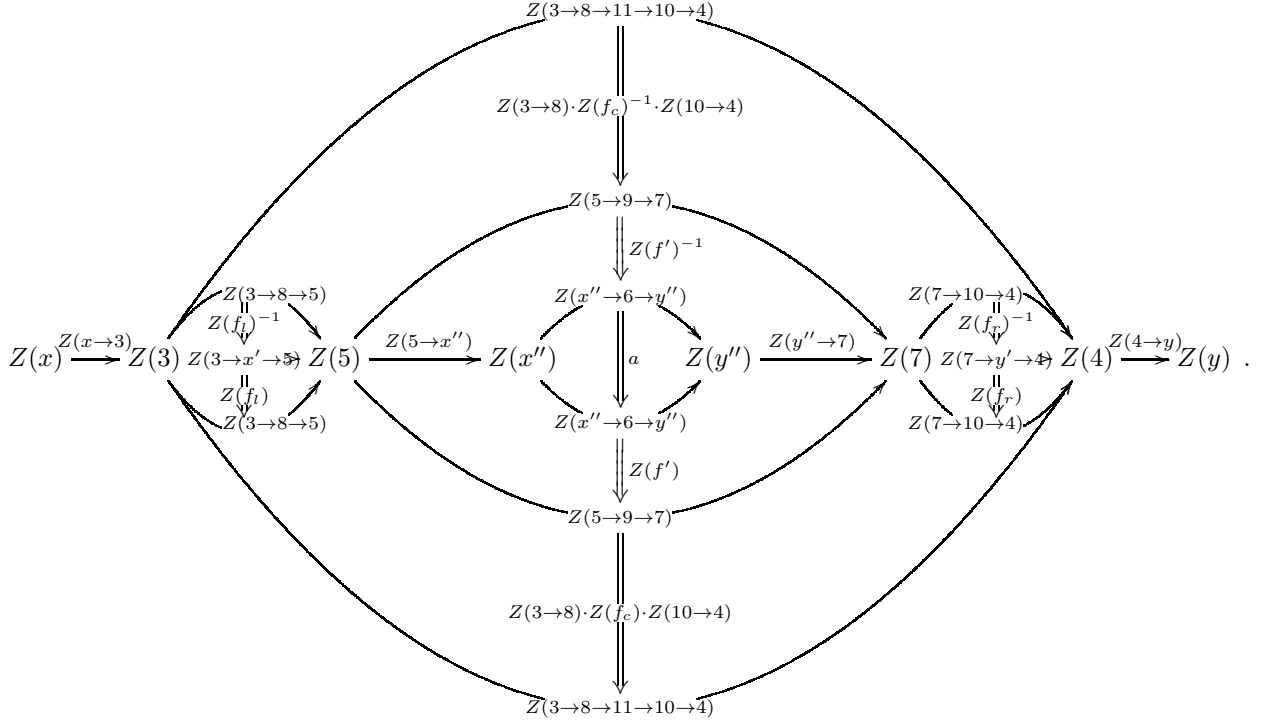
$$
  \xymatrix@C=18pt{
    Z(x)
    \ar[r]^{Z(x \to 3)}
    &
    Z(3)
    \ar@/^2.3pc/[rr]|>>>>>>>{Z(3 \to 8 \to 5)}_{\ }="q1"
    \ar[rr]|{Z(3 \to x' \to 5)}^{\ }="r1"_{\ }="q2"
    \ar@/_2.3pc/[rr]|>>>>>>>{Z(3 \to 8 \to 5)}^{\ }="r2"
    \ar@/^11pc/[rrrrrrrrrr]|{Z(3 \to 8 \to 11 \to 10 \to 4)}_{\ }="v1"
    \ar@/_11pc/[rrrrrrrrrr]|{Z(3 \to 8 \to 11 \to 10 \to 4)}^{\ }="w2"
    &&
    Z(5)
    \ar[rr]^{Z(5 \to x'')}
    \ar@/^5pc/[rrrrrr]|{Z(5  \to 9 \to 7 )}_{\ }="s2"^{\ }="w1"
    \ar@/_5pc/[rrrrrr]|{Z(5  \to 9 \to 7 )}^{\ }="t3"_{\ }="v2"
    &&
    Z(x'')
    \ar@/^2pc/[rr]|{Z(x'' \to 6 \to y'')}_{\ }="s"^{\ }="t2"
    \ar@/_2pc/[rr]|{Z(x'' \to 6 \to y'')}^{\ }="t"_{\ }="s3"
    &&
    Z(y'')
    \ar[rr]^{Z(y'' \to 7)}
    &&
    Z(7)
    \ar@/^2.3pc/[rr]|<<<<<<<{Z(7 \to 10 \to 4)}_{\ }="q3"
    \ar[rr]|{Z(7 \to y' \to 4)}^{\ }="r3"_{\ }="q4"
    \ar@/_2.3pc/[rr]|<<<<<<<{Z(7 \to 10 \to 4)}^{\ }="r4"
    &&
    Z(4)
    \ar[r]^{Z(4 \to y)}
    &
    Z(y)
    \ar@{=>}^a "s"; "t"
    \ar@{=>}^{Z(f')^{-1}} "s2"; "t2"
    \ar@{=>}^{Z(f')} "s3"; "t3"
    \ar@{=>}|{Z(f_l)^{-1}} "q1"; "r1"
    \ar@{=>}|{Z(f_l)} "q2"; "r2"
    \ar@{=>}|{Z(f_r)^{-1}} "q3"; "r3"
    \ar@{=>}|{Z(f_r)} "q4"; "r4"
    \ar@{=>}|{Z(3 \to 8)\cdot Z(f_c)^{-1} \cdot Z(10 \to 4)} "v1"; "w1"
    \ar@{=>}|{Z(3 \to 8)\cdot Z(f_c) \cdot Z(10 \to 4)} "v2"; "w2"
  }
  \,.
$$ The contributions from $f_l$ and $f_r$ manifestly cancel and we are left with
the pasting diagram for the direct inclusion $$
  \mathrm{End}_C(Z(x'' \to 6 \to y''))
  \hookrightarrow
  \mathrm{End}_C(Z(x \to 3 \to 8 \to 11 \to 10 \to 4 \to y))
  \,.
$$ This shows that $$
  \xymatrix{
    \mathcal{A}_Z(O'')
    \ar@{^{(}->}[rr]
    \ar@{^{(}->}[dr]
    &&
    \mathcal{A}_Z(O)
    \\
    &
    \mathcal{A}_Z(O')
    \ar@{^{(}->}[ur]
  }
$$ commutes, as desired.

 To see locality,
 let $O_{x,y}$ and $O_{x',y'}$ be two spacelike separated
  causal subsets inside $O_{(3,5')}$. The relevant
  pasting diagram in $P_2(\mathbb{R}^2)$ is of the form
  $$
  \raisebox{60pt}
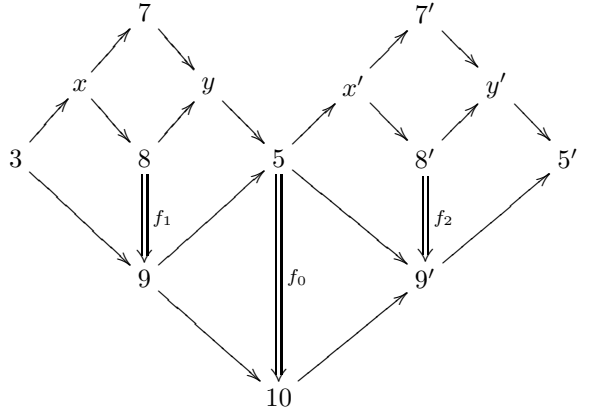
{
\xymatrix@R=13pt@C=13pt{
   &&&
   7
   \ar[dr]
   &&&&
   7'
   \ar[dr]
   \\
   & & x
   \ar[ur]
   \ar[dr]
    && y
   \ar[dr]
   &&
   x'
   \ar[ur]
   \ar[dr]
   && y'
   \ar[dr]
   \\
   & 3
   \ar[ddrr]
   \ar[ur]
   && 8
    \ar[ur]
    \ar@{=>}[dd]^{f_1}
   && 5
   \ar[ddrr]
   \ar@{=>}[dddd]^{f_0}
   \ar[ur]
   &&
   8'
   \ar[ur]
   \ar@{=>}[dd]^{f_2}
   &&
   5'
   \\
   \\
   &&& 9
   \ar[uurr]
   \ar[ddrr]
   &&&&
   9'
   \ar[uurr]
   \\
   \\
   &&&&&
    10
   \ar[uurr]
}
 }
\,. $$ (We are displaying a very symmetric configuration only for notational
convenience. The argument does not depend on that symmetry but just on the fact
that $O_{x,y}$ does not intersect the past of $O_{x',y'}$ and vice versa.) Now
given any two endomorphisms $
  \xymatrix{
    Z(x)
    \ar@/^2pc/[rr]^{Z(x \to 8 \to y)}_{\ }="s"
    \ar@/_2pc/[rr]_{Z(x \to 8 \to y)}^{\ }="t"
    &&
    Z(y)^{a}
    \ar@{=>}^a "s"; "t"
  }
$ and $
  \xymatrix{
    Z(x')
    \ar@/^2pc/[rr]^{Z(x' \to 8' \to y')}_{\ }="s"
    \ar@/_2pc/[rr]_{Z(x' \to 8' \to y')}^{\ }="t"
    &&
    Z(y')
    \ar@{=>}^{a'} "s"; "t"
  }
$ we can either first include $a$ in $\mathrm{End}_C(Z(3 \to 9 \to 10 \to 9' \to
5'))$ and then $a'$, or the other way around. Either way, the total endomorphism
in $\mathrm{End}_C(Z(3 \to 9 \to 10 \to 9' \to 5'))$ is $$
  
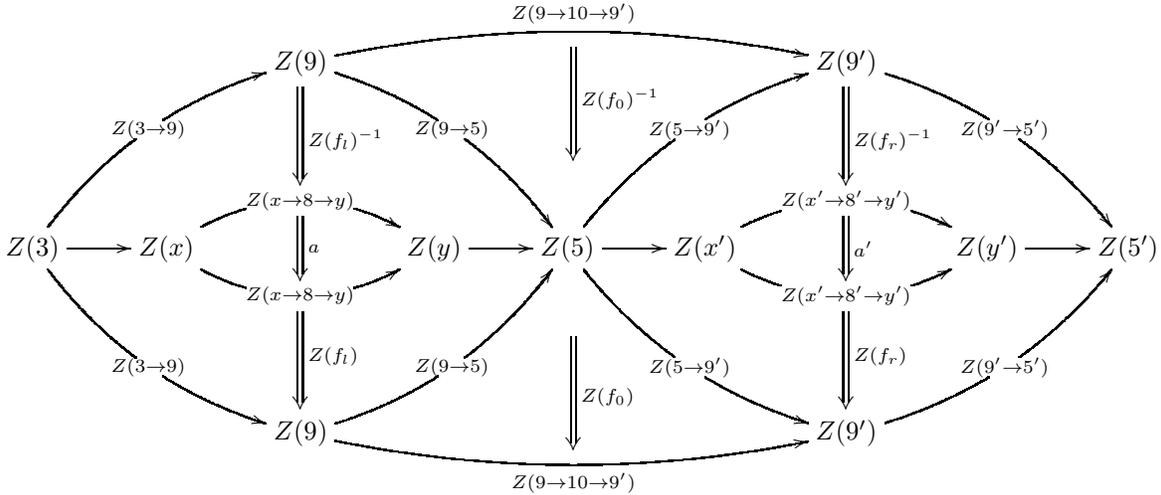
\xymatrix{
    &&
    Z(9)
    \ar@/^1pc/[rrdd]|{Z(9 \to 5)}
    \ar@/^1pc/[rrrr]^{Z(9 \to 10 \to 9')}_{\ }="s9"
    &&&&
    Z(9')
    \ar@/^1pc/[rrdd]|{Z(9' \to 5')}_{\ }="s7"
    \\
    \\
    Z(3)
    \ar[r]
    \ar@/^1pc/[rruu]|{Z(3 \to 9 )}
    \ar@/_1pc/[rrdd]|{Z(3 \to 9 )}
    &
    Z(x)
    \ar@/^1.5pc/[rr]|{Z(x \to 8 \to y)}_{\ }="s1"^{\ }="t3"
    \ar@/_1.5pc/[rr]|{Z(x \to 8 \to y)}^{\ }="t1"_{\ }="s4"
    &
    &
    Z(y)
    \ar[r]
    &
    Z(5)
    \ar[r]
    \ar@/^1pc/[rruu]|{Z(5 \to 9')}
    \ar@/_1pc/[rrdd]|{Z(5 \to 9')}
    &
    Z(x')
    \ar@/^1.5pc/[rr]|{Z(x' \to 8' \to y')}_{\ }="s2"^{\ }="t5"
    \ar@/_1.5pc/[rr]|{Z(x' \to 8' \to y')}^{\ }="t2"_{\ }="s6"
    &
    &
    Z(y')
    \ar[r]
    &
    Z(5')
    \\
    \\
    &&
    Z(9)
    \ar@/_1pc/[rruu]|{Z(9 \to 5)}
    \ar@/_1pc/[rrrr]_{Z(9 \to 10 \to 9')}^{\ }="t8"
    &&&&
    Z(9')
    \ar@/_1pc/[rruu]|{Z(9' \to 5')}
    \ar@{=>}^{a} "s1"; "t1"
    \ar@{=>}^{a'} "s2"; "t2"
    \ar@{=>}^{Z(f_l)^{-1}} "t3"+(0,14); "t3"
    \ar@{=>}^{Z(f_l)} "s4"; "s4"+(0,-14)
    \ar@{=>}^{Z(f_r)^{-1}} "t5"+(0,14); "t5"
    \ar@{=>}^{Z(f_r)} "s6"; "s6"+(0,-14)
    \ar@{=>}^{Z(f_0)^{­-1}} "s9"; "s9"+(0,-16)
    \ar@{=>}^{Z(f_0)} "t8"+(0,16); "t8"
  }
  \,.
$$ This means that the inclusions of $a$ and $a'$ in $\mathrm{End}_C(Z(3 \to 9
\to 10 \to 9' \to 5'))$ commute. \endofproof

\begin{theorem}
  \label{2funct construction}
  This construction extends to a 2-functor
  $$
    \mathbf{FQFT}_{\mathrm{isos}}(\mathbb{R}^2,C)
    \to
    \mathbf{AQFT}(\mathbb{R}^2)
  $$
  faithful on 1-morphisms and trivial on 2-morphisms.
\end{theorem} \proof
  The proof is very analogous to the
  proof of theorem \ref{covariance from equivariance}
  in the next section, only slightly simpler.
\endofproof

\section{Covariance/Equivariance} \label{equivariance}

We had seen definitions for equivariance (``covariance'') of local nets and of
FQFT 2-functors. The following theorem says that these notions are compatible
under our relation of the two.

\begin{theorem}
  \label{covariance from equivariance}
  Every $G$-equivariant structure, definition \ref{equivariance of 2-functors},

  on the FQFT $Z : P_2(\mathbb{R}^2) \to C$
  induces a $G$-equivariant structure, definition \ref{netequi},
  on the AQFT $\mathcal{A}_Z$
  obtained from it according to definition \ref{nets from functors}.
\end{theorem} \proof
  For any $g \in G$ the component map of the pseudonatural transformation
  $f_g$ is
  $$
    f_g
    :
    (\xymatrix{
      x \ar[r]^\gamma & y
    })
    \;
    \mapsto
    \;
    \raisebox{40pt}{
    \xymatrix{
      Z(x)
      \ar[rr]^{Z(\gamma)}
      \ar[dd]_{f_g(x)}^>{\ }="t"
      &&
      Z(y)
      \ar[dd]^{f_g(y)}_<{\ }="s"
      \\
      \\
       Z(g(x))
      \ar[rr]^{Z(g(\gamma))}
      &&
      Z(g(y))
      \ar@{=>}^{f_g(\gamma)} "s"; "t"
   }
   }
   \,.
  $$
  For $\gamma$ the target boundary of the causal subset $O$,
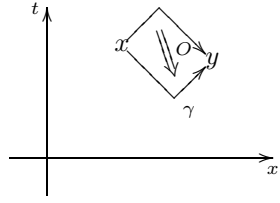
     $$
    \xy
       (10,15)*{x}="s";
       (22,13)*{y}="t";
       (15,20)*{}="top";
       (17,8)*{}="bottom";
       \ar@{<-}_<{t} (0,20); (0,-5)
       \ar@{->}_>{x} (-5,0); (30,0)
       \ar@{-} "s"; "top"
       \ar@{->} "top"; "t"
       \ar@{-} "s"; "bottom"
       \ar@{->}_<\gamma "bottom"; "t"
       \ar@{=>}^O "top"+(0,-3); "bottom"+(0,3)
    \endxy
    $$
  conjugating with the components on the right defines the monoid isomorphism
  $$
    r_g(O) : \mathrm{End}_C(Z(\gamma)) \to \mathrm{End}_C(Z(g(\gamma)))
  $$
  $$
    r_g(O) :
    \left(
       \xymatrix{
         Z(x)
         \ar@/^1pc/[rr]^{Z(\gamma)}_{\ }="s"
         \ar@/_1pc/[rr]_{Z(\gamma)}^{\ }="t"
         &&
         Z(y)
         \ar@{=>}^{a} "s"; "t"
       }
    \right)
    \hspace{6pt}
      \mapsto
    \hspace{6pt}
    \raisebox{70pt}{
    \xymatrix{
       Z(g(x))
      \ar[rr]^{Z(g(\gamma))}
      \ar[dd]|{f_g(x)^{-1}}^>{\ }="t2"
      \ar@/_2pc/[dddd]_{\mathrm{Id}}
      &&
      Z(g(y))
      \ar[dd]|{f_g(y)^{-1}}_<{\ }="s2"
      \ar@/^2pc/[dddd]^{\mathrm{Id}}
      \\
      \\
      Z(x)
      \ar@/^1pc/[rr]|{Z(\gamma)}_{\ }="s1"
      \ar@/_1pc/[rr]|{Z(\gamma)}^{\ }="t1"
      \ar[dd]|{f_g(x)}^>{\ }="t"
      &&
      Z(y)
      \ar[dd]|{f_g(y)}_<{\ }="s"
      \\
      \\
       Z(g(x))
      \ar[rr]^{Z(g(\gamma))}
      &&
      Z(g(y))
      \ar@{=>}^{f_g(\gamma)} "s"; "t"
      \ar@{=>}^{f_g(\gamma)^{-1_p}} "s2"; "t2"
      \ar@{=>}^{a} "s1"; "t1"
   }
   }
   \,.
  $$
  Here $f_g(\gamma)^{-1_p}$ denotes the inverse of the 2-cell $f_g(\gamma)$
  with respect to vertical pasting (which is the ordinary inverse up to
  a re-whiskering).

  We need to check that this construction
  \begin{enumerate}
    \item
      yields a morphism of nets in that
      it makes for all $O' \subset O$ the naturality squares
   $$
    \xymatrix{
      \mathcal{A}_Z(O')
      \ar[r]^{r_g(O')}
      \ar@{^{(}->}[d]
      &
      \mathcal{A}_Z(g(O'))
      \ar@{^{(}->}[d]
      \\
      \mathcal{A}_Z(O)
      \ar[r]^{r_g(O)}
      &
      \mathcal{A}_Z(g(O))
    }
  $$
  commute;
  \item
    produces the commuting triangles in definition \ref{netequi}.
  \end{enumerate}
  This can be seen as follows.
  \begin{enumerate}
  \item

     The pseudo-naturality condition on the components of $f_g$
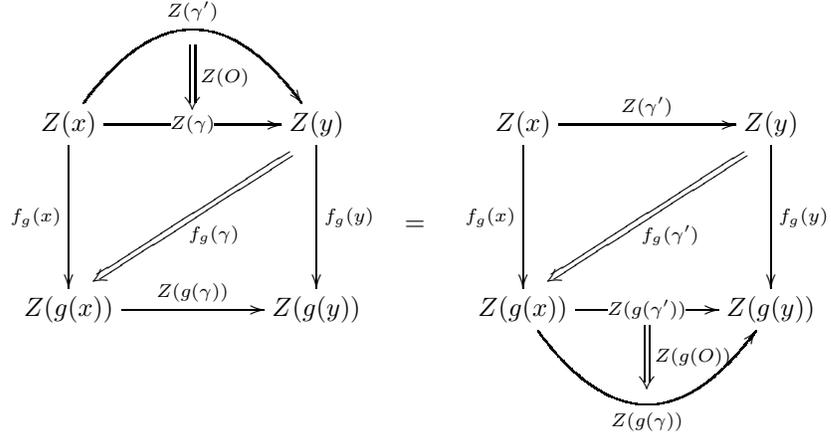
  $$
    \raisebox{40pt}{
    \xymatrix{
      Z(x)
      \ar@/^3pc/[rr]^{Z(\gamma')}_{\ }="s1"
      \ar[rr]|{Z(\gamma)}^{\ }="t1"
      \ar[dd]_{f_g(x)}^>{\ }="t"
      &&
      Z(y)
      \ar[dd]^{f_g(y)}_<{\ }="s"
      \\
      \\
       Z(g(x))
      \ar[rr]^{Z(g(\gamma))}
      &&
      Z(g(y))
      \ar@{=>}^{f_g(\gamma)} "s"; "t"
      \ar@{=>}^{Z(O)} "s1"; "t1"
   }
   }
   \hspace{6pt}
   =
   \hspace{6pt}
    \raisebox{40pt}{
    \xymatrix{
      Z(x)
      \ar[rr]^{Z(\gamma')}^{\ }="t1"
      \ar[dd]_{f_g(x)}^>{\ }="t"
      &&
      Z(y)
      \ar[dd]^{f_g(y)}_<{\ }="s"
      \\
      \\
       Z(g(x))
      \ar[rr]|{Z(g(\gamma'))}_{\ }="s1"
      \ar@/_3pc/[rr]_{Z(g(\gamma))}^{\ }="t1"
      &&
      Z(g(y))
      \ar@{=>}^{f_g(\gamma')} "s"; "t"
      \ar@{=>}^{Z(g(O))} "s1"; "t1"
   }
   }
  $$
  for all $O$
  implies precisely the condition $r_g(O)|_{A(O')} = r_g(O')$
  when applied to our definition \ref{nets from functors} of the inclusion
  map
  $A(O') \hookrightarrow A(O)$: that inclusion was obtained by conjugating
  with
$$
 \raisebox{34pt}
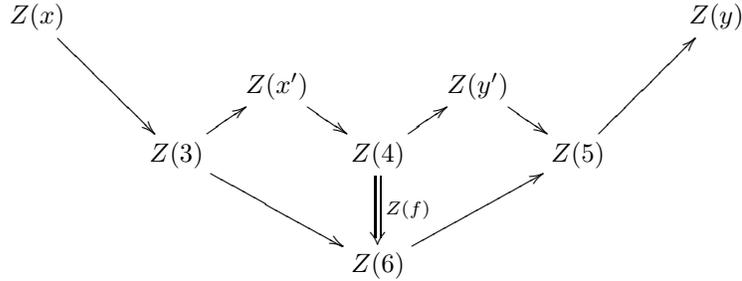
{
\xymatrix@C=10pt@R=10pt{
   Z(x)
   \ar[ddrr]
   &&&&
   &&&&
   Z(y)
   \\
   && & Z(x')  \ar[dr]
    && Z(y')
   \ar[dr]
   &&
   \\
   && Z(3)
   \ar[ddrr]
   \ar[ur]
   && Z(4)
    \ar[ur]
    \ar@{=>}[dd]^{Z(f)}
   && Z(5)
   \ar[uurr]
   \\
   \\
   &&&& Z(6)
   \ar[uurr]
}
 }
  \,.
$$ Following this by the action of $r_g(O)$ amounts to conjugating with $$
 \raisebox{34pt}
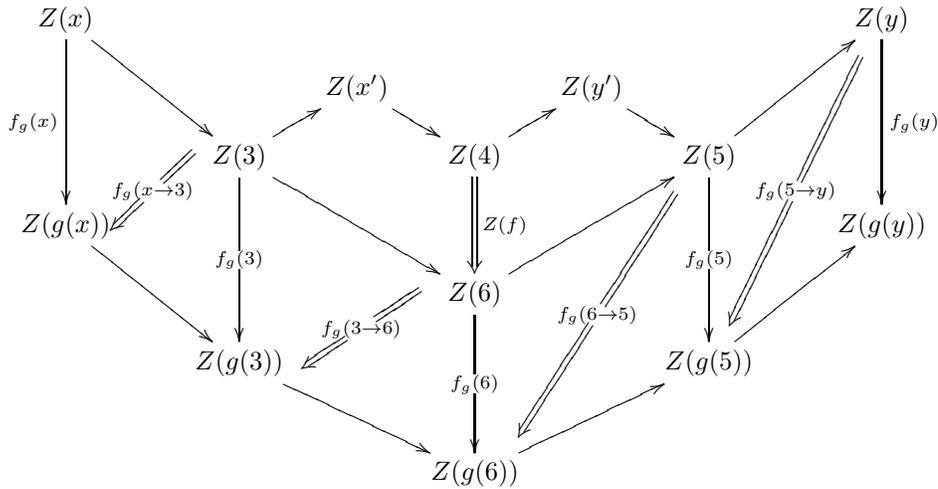
{
\xymatrix@C=10pt@R=10pt{
   Z(x)
   \ar[ddrr]_>{\ }="b"
   \ar[ddd]_{f_g(x)}
   &&&&
   &&&&
   Z(y)
   \ar[ddd]^{f_g(y)}
   \\
   && & Z(x')  \ar[dr]
    && Z(y')
   \ar[dr]
   &&
   \\
   && Z(3)
   \ar[ddrr]_>{\ }="d"
   \ar[ur]
   \ar[ddd]|{f_g(3)}^>{\ }="c"
   && Z(4)
    \ar[ur]
    \ar@{=>}[dd]^{Z(f)}
   && Z(5)
   \ar[uurr]_>{\ }="h"
   \ar[ddd]|{f_g(5)}
   \\
   Z(g(x))
   \ar[ddrr]^<{\ }="a"
   &&&&&&&&
   Z(g(y))
   &&&
   \\
   &&&& Z(6)
   \ar[uurr]_>{\ }="f"
   \ar[ddd]|{f_g(6)}
   \\
   &&Z(g(3))
   \ar[ddrr]^<{\ }="c"
   &&&&
   Z(g(5))
   \ar[uurr]^<{\ }="g"
   &&&&&&
   \\
   &&&&&&&&&&&&
   \\
   &&&& Z(g(6))
   \ar[uurr]^<{\ }="e"
   &&&&&&&&
   \\
   &&&&&&&&&&&&
   \ar@{=>}|{f_g(x \to 3)} "b"; "a"
   \ar@{=>}|{f_g(3\to 6)} "d"; "c"
   \ar@{=>}|{f_g(6 \to 5)} "f"; "e"
   \ar@{=>}|{f_g(5 \to y)} "h"; "g"
 }
 }
  \,.
$$ By pseudonaturality of $f_g$ this equals conjuation with $$
 \raisebox{34pt}
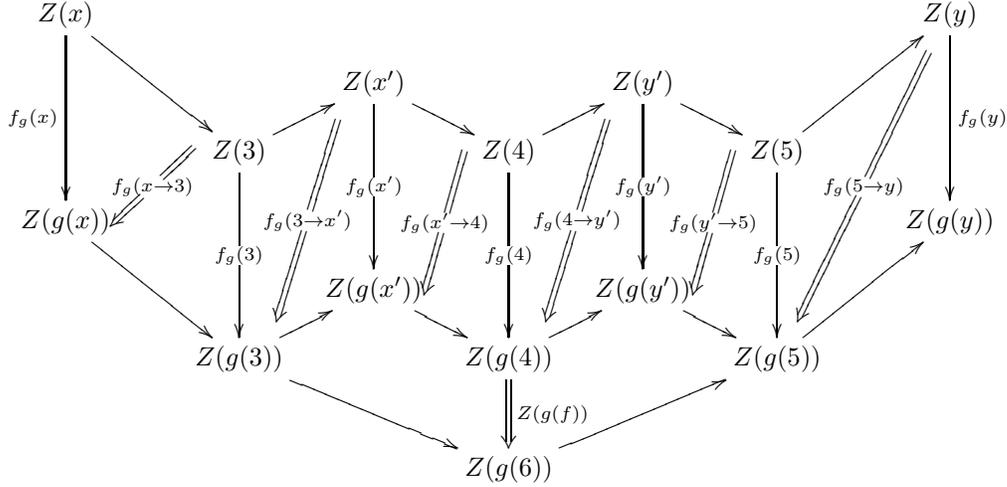
{
\xymatrix@C=10pt@R=10pt{
   Z(x)
   \ar[ddrr]_>{\ }="b"
   \ar[ddd]_{f_g(x)}
   &&&&
   &&&&
   Z(y)
   \ar[ddd]^{f_g(y)}
   \\
   && & Z(x')
   \ar[dr]_>{\ }="r2"
   \ar[ddd]|{f_g(x')}
    && Z(y')
   \ar[dr]_>{\ }="r4"
      \ar[ddd]|{f_g(y')}
   &&
   \\
   && Z(3)
   \ar[ur]_>{\ }="r1"
   \ar[ddd]|{f_g(3)}^>{\ }="c"
   &&
   Z(4)
    \ar[ur]_>{\ }="r3"
   \ar[ddd]|{f_g(4)}
   && Z(5)
   \ar[uurr]_>{\ }="h"
   \ar[ddd]|{f_g(5)}
   \\
   Z(g(x))
   \ar[ddrr]^<{\ }="a"
   &&&&&&&&
   Z(g(y))
   &&&
   \\
   &&&
    Z(g(x'))
    \ar[dr]^<{\ }="s2"
   &
   &
   Z(g(y'))
   \ar[dr]^<{\ }="s4"
   \\
   &&Z(g(3))
   \ar[ddrr]^<{\ }="c"
   \ar[ur]^<{\ }="s1"
   &&
   Z(g(4))
   \ar[ur]^<{\ }="s3"
   \ar@{=>}[dd]^{Z(g(f))}
   &&
   Z(g(5))
   \ar[uurr]^<{\ }="g"
   &&&&&&
   \\
   &&&&&&&&&&&&
   \\
   &&&& Z(g(6))
   \ar[uurr]^<{\ }="e"
   &&&&&&&&
   \\
   &&&&&&&&&&&&
   \ar@{=>}|{f_g(x \to 3)} "b"; "a"
   %\ar@{=>}|{f_g(3\to 6)} "d"; "c"
   %\ar@{=>}|{f_g(6 \to 5)} "f"; "e"
   \ar@{=>}|{f_g(5 \to y)} "h"; "g"
   \ar@{=>}|{f_g(3\to x')} "r1"; "s1"
   \ar@{=>}|{f_g(x'\to 4)} "r2"; "s2"
   \ar@{=>}|{f_g(4\to y')} "r3"; "s3"
   \ar@{=>}|{f_g(y'\to 5)} "r4"; "s4"
 }
 }
  \,.
$$ Since the endomorphism $a$ to be conjugated is localized on $Z(x') \to
Z(y')$ 
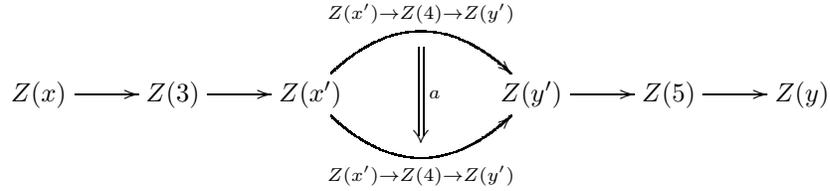
$$
  \xymatrix{
    Z(x)
    \ar[r]
    &
    Z(3)
    \ar[r]
    &
    Z(x')
    \ar@/^2pc/[rr]^{Z(x') \to Z(4) \to Z(y')}_{\ }="s"
    \ar@/_2pc/[rr]_{Z(x') \to Z(4) \to Z(y')}^{\ }="t"
    &&
    Z(y')
    \ar[r]
    &
    Z(5)
    \ar[r]
    &
    Z(y)
    \ar@{=>}^a "s"; "t"
  }
$$ both $f_g(x \to 3 \to x')$ and $f_g(y'\to 5 \to y)$ drop out when
conjugating and only conjugation with $f_g(x \to 4 \to y')$ acts
nontrivially. But that precisely amounts to first applying $r_g(O')$ and
then injecting into $O$.

  \item
  The equivariance triangle condition in definition
  \ref{equivariance of 2-functors}
  says precisely that $r_g(O)$ makes the required
  covariance triangle in definition \ref{netequi}
  commute: To see this it is convenient to
  equivalently rewrite the previous equation for $r_g(O)$ as
  $$
    \raisebox{40pt}
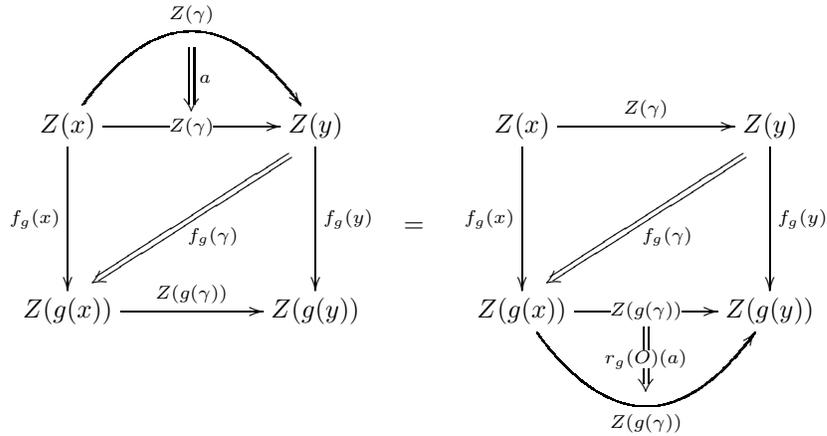
{
    \xymatrix{
      Z(x)
      \ar@/^3pc/[rr]^{Z(\gamma)}_{\ }="s1"
      \ar[rr]|{Z(\gamma)}^{\ }="t1"
      \ar[dd]_{f_g(x)}^>{\ }="t"
      &&
      Z(y)
      \ar[dd]^{f_g(y)}_<{\ }="s"
      \\
      \\
      Z(g(x))
      \ar[rr]^{Z(g(\gamma))}
      &&
      Z(g(y))
      \ar@{=>}^{f_g(\gamma)} "s"; "t"
      \ar@{=>}^{a} "s1"; "t1"
   }
   }
   \hspace{6pt}
   =
   \hspace{6pt}
    \raisebox{40pt}{
    \xymatrix{
      Z(x)
      \ar[rr]^{Z(\gamma)}^{\ }="t1"
      \ar[dd]_{f_g(x)}^>{\ }="t"
      &&
      Z(y)
      \ar[dd]^{f_g(y)}_<{\ }="s"
      \\
      \\
       Z(g(x))
      \ar[rr]|{Z(g(\gamma))}_{\ }="s1"
      \ar@/_3pc/[rr]_{Z(g(\gamma))}^{\ }="t1"
      &&
      Z(g(y))
      \ar@{=>}^{f_g(\gamma)} "s"; "t"
      \ar@{=>}|{r_g(O)(a)} "s1"; "t1"
   }
   }
  $$
  for all $a \in \mathrm{End}(Z(\gamma))$.
  Accordingly, we have for the composition of two transformations
  $$
    \raisebox{60pt}{
    \xymatrix{
      Z(x)
      \ar@/^3pc/[rr]^{Z(\gamma)}_{\ }="s1"
      \ar[rr]|{Z(\gamma)}^{\ }="t1"
      \ar[dd]_{f_{g_1}(x)}^>{\ }="t"
      &&
      Z(y)
      \ar[dd]^{f_{g_1}(y)}_<{\ }="s"
      \\
      \\
      Z({g_1}(x))
      \ar[rr]^{Z({g_1}(\gamma))}
      \ar[dd]_{f_{g_2}(g_1(x))}^>{\ }="t2"
      &&
      Z(g_1(y))
      \ar[dd]^{f_{g_2}(g_1(y))}_<{\ }="s2"
      \\
      \\
      Z({(g_1 g_2)}(x))
      \ar[rr]^{Z({g_1g_2}(\gamma))}
      &&
      Z((g_1g_2)(y))
      \ar@{=>}^{f_{g_1}(\gamma)} "s"; "t"
      \ar@{=>}^{a} "s1"; "t1"
      \ar@{=>}^{f_{g_2}(g_1(\gamma))} "s2"; "t2"
   }
   }
   \hspace{6pt}
   =
   \hspace{6pt}
    \raisebox{60pt}{
    \xymatrix{
      Z(x)
      \ar[rr]^{Z(\gamma)}
      \ar[dd]_{f_{g_1}(x)}^>{\ }="t"
      &&
      Z(y)
      \ar[dd]^{f_{g_1}(y)}_<{\ }="s"
      \\
      \\
      Z(x)
      \ar[rr]^{Z(g_1(\gamma))}
      \ar[dd]_{f_{g_2}(g_1(x))}^>{\ }="t2"
      &&
      Z(y)
      \ar[dd]^{f_{g_2}(g_1(y))}_<{\ }="s2"
      \\
      \\
       Z({(g_1g_2)}(x))
      \ar[rr]|{Z({g_1g_2}(\gamma))}_{\ }="s1"
      \ar@/_3pc/[rr]_{Z({g_1 g_2}(\gamma))}^{\ }="t1"
      &&
      Z({(g_1 g_2)}(y))
      \ar@{=>}^{f_{g_1}(\gamma)} "s"; "t"
      \ar@{=>}^{f_{g_2}(g_1(\gamma))} "s2"; "t2"
      \ar@{=>}|{r_{g_2}(g_1(O))\circ r_{g_1}(O)(a)} "s1"; "t1"
   }
   }
  $$
  for all $a \in \mathrm{End}(Z(\gamma))$. Using now the triangle
  of pseudonatural transformations in definition \ref{equivariance of
  2-functors}
  this is equivalent to
  $$
    \hspace{-6cm}
    \raisebox{60pt}{
    \xymatrix{
      &
      Z(x)
      \ar@/^3pc/[rr]|{Z(\gamma)}_{\ }="s1"
      \ar[rr]|{Z(\gamma)}^{\ }="t1"
      \ar[ddl]|{f_{g_1}(x)}^>{\ }
      \ar[dddd]|{f_{g_1 g_2}(x)}^>{\ }="tu"_{\ }="sv"
      &&
      Z(y)
      \ar[ddr]|{f_{g_1}(y)}_{\ }="s"
      \ar[dddd]|{f_{g_1 g_2}(y)}_<{\ }="su"^{\ }="t"
      \\
      \\
      Z(g_1(x))
      \ar[ddr]|{f_{g_2}(g_1(x))}^{\ }="ta"
      &&
      &&Z(g_1(y))
      \ar[ddl]^{f_{g_2}(g_1(y))}_<{\ }="sa"
      \\
      \\
      &
      Z((g_1 g_2)(x))
      \ar[rr]^{Z_{(g_1 g_2)(\gamma)}}
      &&
      Z((g_1 g_2)(y))
      \ar@{=>}|{F_{g_1,g_2}(y)} "s"; "t"
      %\ar@{=>}^{f_{g_1}(\gamma)^{-1}} "s2"; "t2"
      \ar@{=>}|{f_{g_1 g_2}(\gamma)} "su"; "tu"
      %\ar@{=>}|{f_{g_1 g_2}(\gamma)^{-1}} "sb"; "tb"
      \ar@{=>}^{a} "s1"; "t1"
      \ar@{=>}|{F_{g_1,g_2}(x)^{-1}} "sv"; "ta"
      %\ar@{=>} "sl"; "tl"
      %\ar@{=>} "sn"; "tn"
    }
    }
    \hspace{6pt}
     =
  $$
  $$
    \hspace{6cm}
      \raisebox{60pt}{
      \xymatrix{
      &
      Z(x)
      \ar[rr]|{Z(\gamma)}
      \ar[ddl]|{f_{g_1}(x)}^>{\ }
      \ar[dddd]|{f_{g_1 g_2}(x)}^>{\ }="tu"_{\ }="sv"
      &&
      Z(y)
      \ar[ddr]|{f_{g_1}(y)}_{\ }="s"
      \ar[dddd]|{f_{g_1 g_2}(y)}_<{\ }="su"^{\ }="t"
      \\
      \\
      Z(g_1(x))
      \ar[ddr]|{f_{g_2}(g_1(x))}^{\ }="ta"
      &&
      &&Z(g_1(y))
      \ar[ddl]^{f_{g_2}(g_1(y))}_<{\ }="sa"
      \\
      \\
      &
      Z((g_1 g_2)(x))
      \ar[rr]|{Z_{(g_1 g_2)(\gamma)}}_{\ }="s1"
      \ar@/_3pc/[rr]|{Z_{(g_1 g_2)(\gamma)}}_{\ }="t1"
      &&
      Z((g_1 g_2)(y))
      \ar@{=>}|{F_{g_1,g_2}(y)} "s"; "t"
      %\ar@{=>}^{f_{g_1}(\gamma)^{-1}} "s2"; "t2"
      \ar@{=>}|{f_{g_1 g_2}(\gamma)} "su"; "tu"
      %\ar@{=>}|{f_{g_1 g_2}(\gamma)^{-1}} "sb"; "tb"
      \ar@{=>}|{r_{g_2}(g_1(O))\circ g_{g_1}(O)a} "s1"; "t1"
      \ar@{=>}|{F_{g_1,g_2}(x)^{-1}} "sv"; "ta"
      %\ar@{=>} "sl"; "tl"
      %\ar@{=>} "sn"; "tn"
    }
    }
    \,.
  $$
  But in this equation we can cancel the $F_{\cdot,\cdot}$ on both sides
  to obtain
  $$
    \raisebox{60pt}{
    \xymatrix{
      Z(x)
      \ar@/^3pc/[rr]|{Z(\gamma)}_{\ }="s1"
      \ar[rr]|{Z(\gamma)}^{\ }="t1"
      \ar[dddd]|{f_{g_1 g_2}(x)}^>{\ }="tu"_{\ }="sv"
      &&
      Z(y)
      \ar[dddd]|{f_{g_1 g_2}(y)}_<{\ }="su"^{\ }="t"
      \\
      \\
      \\
      \\
      Z((g_1 g_2)(x))
      \ar[rr]^{Z_{(g_1 g_2)(\gamma)}}
      &&
      Z((g_1 g_2)(y))
      %
      %\ar@{=>}^{f_{g_1}(\gamma)^{-1}} "s2"; "t2"
      \ar@{=>}|{f_{g_1 g_2}(\gamma)} "su"; "tu"
      %\ar@{=>}|{f_{g_1 g_2}(\gamma)^{-1}} "sb"; "tb"
      \ar@{=>}^{a} "s1"; "t1"
      %\ar@{=>} "sl"; "tl"
      %\ar@{=>} "sn"; "tn"
    }
    }
    \hspace{6pt}
     =
    \hspace{6pt}
      \raisebox{60pt}{
      \xymatrix{
      Z(x)
      \ar[rr]|{Z(\gamma)}
      \ar[dddd]|{f_{g_1 g_2}(x)}^>{\ }="tu"_{\ }="sv"
      &&
      Z(y)
      \ar[dddd]|{f_{g_1 g_2}(y)}_<{\ }="su"^{\ }="t"
      \\
      \\
      \\
      \\
      Z((g_1 g_2)(x))
      \ar[rr]|{Z_{(g_1 g_2)(\gamma)}}_{\ }="s1"
      \ar@/_3pc/[rr]|{Z_{(g_1 g_2)(\gamma)}}_{\ }="t1"
      &&
      Z((g_1 g_2)(y))
      %
      %\ar@{=>}^{f_{g_1}(\gamma)^{-1}} "s2"; "t2"
      \ar@{=>}|{f_{g_1 g_2}(\gamma)} "su"; "tu"
      %\ar@{=>}|{f_{g_1 g_2}(\gamma)^{-1}} "sb"; "tb"
      \ar@{=>}|{r_{g_2}(g_1(O))\circ g_{g_1}(O)a} "s1"; "t1"
      %\ar@{=>} "sl"; "tl"
      %\ar@{=>} "sn"; "tn"
    }
    }
    \,.
  $$
  This shows that $r_{g_2}(g_1(O))\circ r_{g_1}(O)(a) =
  r_{g_1 g_2}(O)(a)$.

  \end{enumerate}

\endofproof

\section{Examples} \label{examples}

\subsection{1-dimensional case} Before looking at concrete examples for 2-FQFTs
on Minkowski space it is again helpful to first recall some simple facts in the
1-dimensional case from our perspective.

We can regard ordinary quantum mechanics as given by an associated $U(E)$-bundle
with connection on the real line (the ``worldline'') for $E$ some Hilbert space.
This bundle is necessarily trivializable. After picking a trivialization its
globally defined $\mathrm{Lie}(U(E))$-valued connection 1-form is $$
  A = i H dt \in \Omega^1(\mathbb{R}^1, \mathfrak{u}(E))
$$ with $t$ the canonical coordinate and $H$ a self-adjoint operator on $E$: the
Hamilton operator. The quantum time evolution operator $$
  Z :
  (\xymatrix{
    t_0 \ar[r] & t_1
  })
  \mapsto
  (\xymatrix{
    E \ar[rr]^{P \exp(\int_{[t_0,t_1]}A )} && E
  })
$$ is nothing but the parallel transport with respect to $A$ (see for instance
\cite{SWI}).

In general $H$ depends on $t$, in which case one speaks of \emph{time dependent}
quantum mechanics and the above formula, with its ``path ordered exponential''
on the right, is what is usually referred to as the \emph{Dyson formula} in
quantum mechanics textbooks. In that case there is no translational invariance
on the worldline.

If however $H$ is constant we have \emph{time independent} quantum mechanics. In
that case the quantum time evolution propagator reads $$
  Z :
  (\xymatrix{
    t_0 \ar[r] & t_1
  })
  \mapsto
  (\xymatrix{
    E \ar[rr]^{P \exp(\int_{[t_0,t_1]}A )} && E
  })
  =
  (\xymatrix{
    E \ar[rr]^{\exp( i (t_1-t_0) H )} && E
  })
  \,.
$$ In either case, there is a canonical equivariant structure, definition
\ref{equivariance of 2-functors}, on $Z$ with respect to the action of
$\mathbb{R}$ on $\mathbb{R}$ by translations: for $a \in \mathbb{R}$ the
components of the natural transformation $$
  \xymatrix{
    Z \ar[r]^{f_t} & a^* Z
  }
$$ are simply $$
  f_a : x \mapsto
  (\xymatrix{
     E_x
     \ar[rr]^{Z(x\to x+a)}
     &&
     E_{x+a}
  })
  \,.
$$ Naturality of $f_t$ and commutativity of the equivariance coherence triangle
both follow directly from the functoriality of $Z$. The equivariant structure on
the net $\mathcal{A}_Z$ induced by this according to section \ref{equivariance}
is that which acts on each local algebra $\mathcal{A}_Z(O_x)$ by the Heisenberg
propagation rule $a \mapsto Z(x \to x+a)\circ a \circ Z(x \to x+a)^{-1}$.

\subsection{Examples from parallel 2-transport}

The above shows that the dynamics of quantum mechanics (1+0-dimensional QFT) can
be entirely thought of as a vector bundle (or Hilbert bundle, rather) with
connection on the ``worldline'' $\mathbb{R}$.

Similarly, 2-vector 2-bundles \cite{Bartels,Wockel} ($\simeq$ gerbes) with
connection \cite{BS,SWII,SWIII,ndclecture} on the ``worldsheet'' $\mathbb{R}^2$
can be regarded as giving the dynamics of (1+1)-dimensional QFT. Indeed, every
parallel transport 2-functor on $\mathbb{R}^2$ as in \cite{BS,SWII,SWIII} gives
an example of a 2-FQFT in the sense definition \ref{eFQFT}, simply by
restricting it from all 2-paths in $\mathbb{R}^2$ to those contained in
$P_2(\mathbb{R}^2)$. From each such 2-functor one obtains, by theorem
\ref{consistency}, a local net of monoids. Whether this local net of monoids has
any covariance depends, according to proposition \ref{covariance from
equivariance}, or whether or not the 2-functor has any equivariant structure.
Whether the net of \emph{monoids} obtained from the 2-functor is actually a net
of algebras with certain extra structure (in particular $C^*$, von Neumann)
depends on what precisely the 2-functor takes values in over 1-morphisms,
because that determines what the endomorphism monoids are like.

While not every 2-bundle on 2-dimensional base space is necessarily
trivializable, we here want to restrict attention to the case that the 2-bundle
is trivializable. (If not, global effects such as described in \cite{FS} will
play a role, too.) Then we can assume its parallel transport 2-functor to come
from globally defined differential form data. If we require the 2-functor to be
\emph{strict} and to take values in a 2-groupoid with a single object, which we
shall denote $\mathbf{B}G$, then theorem 2.20 in \cite{SWII} says that it comes
precisely from a pair consisting of a 1-form and a 2-form $$
  A \in \Omega^1(\mathbb{R}^2,\gg),\;
  B \in \Omega^2(\mathbb{R}^2, \hh)
$$ with values in Lie algebras $\gg$ and $\hh$ which form a differential crossed
module $(\xymatrix{
  \hh \ar[r]^t & \gg \ar[r]^\alpha & \mathrm{der}(\gg)
})$ such that $$
  F_A + t_* \circ B = 0
  \,,
$$ where $F_A \in \Omega^2(\mathbb{R}^2,\gg)$ is the curvature 2-form of $A$. We
write $$
  Z_{(A,B)} : P_2(\mathbb{R}^2) \to \mathbf{B}G
$$ for the 2-functor obtained this way. The local net $
  A_{Z_{(A,B)}}
$ obtained from this by theorem \ref{consistency} is a \emph{local net of
groups}.

We get proper nets of local \emph{algebras} by passing instead to an
\emph{associated} parallel 2-transport functor \cite{SWIII}, 
which is induced by a
2-representation of $G$ on 2-vector space, i.e. a 2-functor $$
  \rho : \mathbf{B}G \to 2\mathrm{Vect}
  \,,
$$ where $2\mathrm{Vect}$ denotes a 2-category of 2-vector spaces. In
particular, \cite{canrep}, there are large classes of 2-representations which
factor through the bicategory of bimodules $$
  \xymatrix{
     \mathbf{B}G
       \ar[dr]
       \ar[rr]^\rho && 2\mathrm{Vect}
     \\
     &
     \mathrm{Bimod}
     \ar[ur]
  }
  \,,
$$ More details on this are summarized in appendix \ref{2veccanrep}
and in \cite{SWIII}.

The corresponding associated 2-FQFT functor $$
  Z_{\rho(A,B)}
  :
  \xymatrix{
     P_2(\mathbb{R}^2)
     \ar[r]^{Z_{(A,B)}}
     &
     \mathbf{B}G
     \ar[r]
     &
     \mathrm{Bimod}
     \ar[r]
     &
     2\mathrm{Vect}
  }
$$ sends each edge to a bimodule over some algebra. 2-Functors of this form and
interpreted as 2-FQFTs have in particular been considered in \cite{ST}.

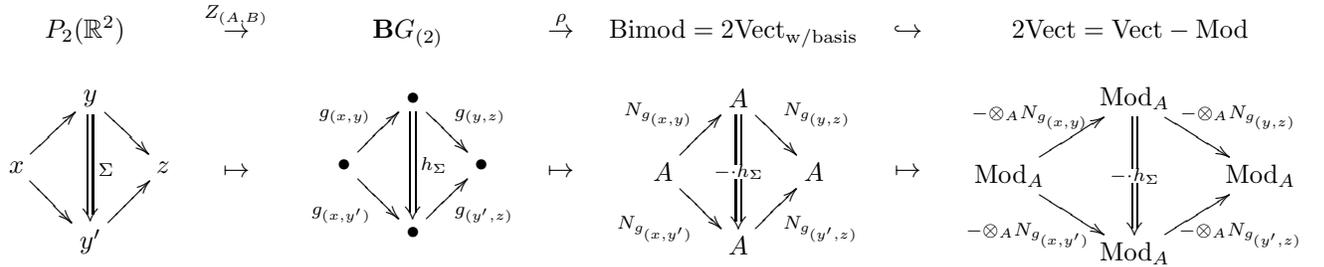
\begin{figure}[h]
  \hspace{-1cm}
  $$
   \begin{array}{ccccccccc}
    P_2(\mathbb{R}^2)
    &
      \stackrel{Z_{(A,B)}}{\to}
    &
    \mathbf{B}G_{(2)}
    &
      \stackrel{\rho}{\to}
    &
    \mathrm{Bimod} = 2\mathrm{Vect}_{\mathrm{w/basis}}
    &
      \stackrel{}{\hookrightarrow}
    &
    2\mathrm{\mathrm{Vect}}
    =
    \mathrm{Vect}-\mathrm{Mod}
    \\
    \\
    \raisebox{30pt}{
    \xymatrix@R=15pt@C=15pt{
      & y
      \ar[dr]
      \ar@{=>}[dd]^\Sigma
      \\
      x
      \ar[ur]
      \ar[dr]
      &&
      z
      \\
      &
      y'
      \ar[ur]
    }
    }
     &
     \mapsto
     &
    \raisebox{30pt}{
    \xymatrix@R=15pt@C=15pt{
      & \bullet
      \ar[dr]^{g_{(y,z)}}
      \ar@{=>}[dd]^{h_\Sigma}
      \\
      \bullet
      \ar[ur]^{g_{(x,y)}}
      \ar[dr]_{g_{(x,y')}}
      &&
      \bullet
      \\
      &
      \bullet
      \ar[ur]_{g_{(y',z)}}
    }
    }
     &
     \mapsto
     &
    \raisebox{30pt}{
    \xymatrix@R=15pt@C=15pt{
      & A
      \ar[dr]^{N_{g_{(y,z)}}}
      \ar@{=>}[dd]|{­-\cdot h_\Sigma}
      \\
      A
      \ar[ur]^{N_{g_{(x,y)}}}
      \ar[dr]_{N_{g_{(x,y')}}}
      &&
      A
      \\
      &
      A
      \ar[ur]_{N_{g_{(y',z)}}}
    }
    }
     &
     \mapsto
     &
    \raisebox{30pt}{
    \xymatrix@R=15pt@C=15pt{
      & \mathrm{Mod}_A
      \ar[dr]^{- \otimes_A N_{g_{(y,z)}}}
      \ar@{=>}[dd]|{­-\cdot h_\Sigma}
      \\
      \mathrm{Mod}_A
      \ar[ur]^{-\otimes_A N_{g_{(x,y)}}}
      \ar[dr]_{- \otimes_A N_{g_{(x,y')}}}
      &&
      \mathrm{Mod}_A
      \\
      &
      \mathrm{Mod}_A
      \ar[ur]_{- \otimes_A N_{g_{(y',z)}}}
    }
    }
   \end{array}
  $$
  \caption{
    2-Vector transport coming from a 2-connection
    $(A,B) \in \Omega^\bullet(\mathbb{R}^2,(\hh \to \gg))$
    with values in the strict Lie 2-algebra $(\hh \to \gg)$
    and the canonical representation $\rho$ of the corresponding
    strict Lie 2-group $G_{(2)}$ on 2-vector spaces.
    The 2-FQFT obtained this way assigns algebras to points,
    bimodules to paths and bimodule homomorphisms to surfaces.
    The corresponding local net $A_{Z_{(A,B)}}$ assigns algebras
    of bimodule endomorphisms.
  }
\end{figure}

Therefore the corresponding local net $A_{Z_{\rho(A,B)}}$ sends each $O \in
S(\mathbb{R}^2)$ to an algebra of bimodule endomorphisms. This is reminiscent of
various other constructions that have been considered in the context of AQFT.
But a more detailed discussion will have to be given elsewhere.

As in the 1-dimensional case, we canonically have an equivariant structure on
$Z$ and on $\mathcal{A}_{Z}$ with respect to any 1-parameter group of
translations which resepcts the light-cone structure. Let in particular
$\mathbb{R}$ act by translation along the canonical time coordinare on
$\mathbb{R}^2$. Then for $a \in \mathbb{R}$ the component of the pseudonatural
transformation $$
  
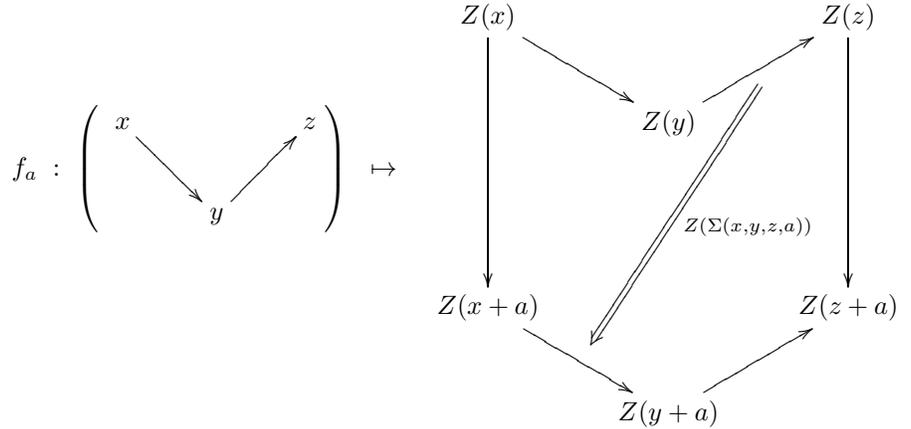
\xymatrix{
    Z
    \ar[rr]^{f_a}
    &&
    a^* Z
  }
$$ is $$
 f_a
 \hspace{3pt}
  :
 \hspace{3pt}
 \left(
   \raisebox{20pt}{
   \xymatrix{
    x \ar[dr]&& z
    \\
    & y
    \ar[ur]
 }}
 \right)
 \hspace{6pt}
   \mapsto
 \hspace{6pt}
 \raisebox{60pt}{
 \xymatrix{
    Z(x) \ar[dr]
    \ar[ddd]
    &&
    Z(z)
    \ar[ddd]
    \\
    & Z(y)
    \ar[ur]_{\ }="s"
    \\
    \\
    Z(x+a)
      \ar[dr]^{\ }="t"
     &&
    Z(z+a)
    \\
    &
    Z(y+a)
    \ar[ur]
    \ar@{=>}^{Z(\Sigma(x,y,z,a))} "s"; "t"
 }
 }
 \,,
$$ where $\Sigma(x,y,z,a)$ denotes the surface swept out by the path $x \to y
\to z$ when translating it continuously to $(x+a) \to (y+a) \to (z+a)$. This
surface is not part of $P_2(\mathbb{R}^2)$ the way we have defined it, but is a
more general 2-path in $\mathbb{R}^2$ on which we can evaluate our 2-functor
$Z$, by assumption.

Pseudonaturality and coherence of the assignment $f_a$ for all $a \in
\mathbb{R}$ is a direct consequence of the 2-functoriality of $Z$, very similar
to the 1-dimensional case. The induced equivariant structure on the net
$\mathcal{A}_Z$ is the local Heisenberg picture time propagation.

\subsection{2-Functors constant on one object}

A simple class of examples worth looking at to get a feeling for the situation
are those FQFT 2-functors $Z$ on $P_2(\mathbb{R}^2)$ which assign a fixed object
$V \in \mathrm{Obj}(C)$ to each point of $\mathbb{R}^2$, send all paths to the
identity morphism on that object and all surfaces to the identity 2-morphism on
this identity 1-morphism.

The local net $\mathcal{A}_Z$ obtained from such a 2-functor is constant. It
assigns the same monoid to all causal subsets: $$
  \mathcal{A}_Z : O \mapsto \mathrm{End}(\mathrm{Id}_V)
  \,.
$$ For this to be a local net, it must be true that
$\mathrm{End}(\mathrm{Id}_V)$ is a commutative monoid. And indeed it is: this is
the Eckmann-Hilton argument which holds in general for 2-endomorphisms of
identity 1-functors. The argument is entirely analogous (and that is of course
no coincidence) to that which shows that the second homotopy group of any space
is abelian.

In \cite{GanterKapranov} the endomorphisms of the identity on an object $V$ in a
2-category $C$ is interpreted as the \emph{trace} of the identity on $V$, which
in turn is interpreted in \cite{BW} as the \emph{dimension} of $V$: $$
  A_Z(O) = \mathrm{End}(\mathrm{Id}_V) =: \mathrm{Tr}(\mathrm{Id}_V)
  =: \mathrm{dim}(V)
  \,.
$$ For instance (see \cite{BW}) if $V = \mathrm{Rep}(H)$ is the category of
representations of some group or groupoid $H$, regarded as a 2-vector space,
then $\mathrm{dim}(V) = Z(\mathbb{C}(H))$ is the center of the group ring of
$H$.

Another example, \cite{GanterKapranov}: if $C$ is the bicategory of bimodules,
$C = \mathrm{Bimod}$, and $V$ is any algebra, then $\mathrm{dim}(V)$ is the 0th
Hochschild cohomology of $V$. Full Hochschild cohomology is obtained by taking
the derived category of bimodules.

Of particular interest are objects $V$ with a representation (meaning:
2-representation!) of the Poincar{\'e} group $G$ in two dimensions, or some
related group, on them. 2-Representations of the Poincar{\'e} group have been
examined for instance in \cite{Poin2rep}. The constant FQFT 2-functor on such an
object canonically carries a nontrivial $G$-equivariant structure in the sense
of section \ref{equivariance}, hence induces a covariant structure on the
corresponding local net.

\subsection{Lattice models}

All our definitions and constructions make sense for $S(\mathbb{R}^2)$  and
$P_2(\mathbb{R}^2)$ replaced by their restrictions $S(\mathbb{Z}^2)$ and
$P_2(\mathbb{Z}^2)$ along that embedding $\mathbb{Z}^2 \hookrightarrow
\mathbb{R}^2$ which makes addition of $(1,0)$ a lightlike translation. This
allows to see a class of important examples without the need to worry about weak
2-categories and issues in functional analysis.

Let $$
  C := \mathbf{B}\mathrm{Vect}
  =
  \left\lbrace
    \xymatrix{
      \bullet
      \ar@/^1.7pc/[rr]^{V}_{\ }="s"
      \ar@/_1.7pc/[rr]_{W}^{\ }="t"
      &&
      \bullet
      \ar@{=>}^\phi "s"; "t"
    }
    \;
    |
    (\xymatrix{V \ar[r]^\phi & W}) \in \mathrm{Vect}
  \right\rbrace
$$ be the strict 2-category obtained from the strict monoidal category of
finite-dimensional vector spaces: it has a single object, its 1-morphisms are
finite dimensional vector spaces with composition of morphisms being the tensor
product of vector spaces, and 2-morphisms are linear maps $\xymatrix{V
\ar[r]^\phi & W}$ between vector spaces.

Pick a fixed finite dimensional vector space $V$ and consider the two 2-FQFT
2-functors $$
  Z_\Vert : P_2(\mathbb{Z}^2) \to \mathbf{B}\mathrm{Vect}
$$ and $$
  Z_\times : P_2(\mathbb{Z}^2) \to \mathbf{B}\mathrm{Vect}
$$ which assign $V$ to every elementary 1-morphism in $P_2(\mathbb{Z}^2)$ and
which assign to every elementary square the linear map $$
  Z_\Vert
  \left(
    \raisebox{40pt}{
    \xymatrix{
       & y
       \ar[dr]
       \ar@{=>}[dd]
       \\
       x
       \ar[ur]
       \ar[dr]
       &&
       z
       \\
       &
       y'
       \ar[ur]
    }
    }
  \right)
  \hspace{7pt}
    :=
  \hspace{7pt}
    \raisebox{40pt}{
    \xymatrix{
       & \bullet
       \ar[dr]^V_{\ }="s2"
       \\
       \bullet
       \ar[ur]^V_{\ }="s1"
       \ar[dr]_V^{\ }="t1"
       &&
       \bullet
       \\
       &
       \bullet
       \ar[ur]_V^{\ }="t2"
       \ar@{=>}|{\mathrm{Id}} "s1"; "t1"
       \ar@{=>}|{\mathrm{Id}} "s2"; "t2"
    }
    }
  \hspace{5pt}
  =
  \hspace{5pt}
    \xymatrix{
      \bullet
      \ar@/^1.7pc/[rr]^{V \otimes V}_{\ }="s"
      \ar@/_1.7pc/[rr]_{V \otimes V}^{\ }="t"
      &&
      \bullet
      \ar@{=>}^{\mathrm{Id}} "s"; "t"
    }
$$ and $$
  Z_\times
  \left(
    \raisebox{40pt}{
    \xymatrix{
       & y
       \ar[dr]
       \ar@{=>}[dd]
       \\
       x
       \ar[ur]
       \ar[dr]
       &&
       z
       \\
       &
       y'
       \ar[ur]
    }
    }
  \right)
  \hspace{7pt}
    :=
  \hspace{7pt}
    \raisebox{40pt}{
    \xymatrix{
       & \bullet
       \ar[dr]^V_{\ }="s2"
       \\
       \bullet
       \ar[ur]^V_{\ }="s1"
       \ar[dr]_V^{\ }="t1"
       &&
       \bullet
       \\
       &
       \bullet
       \ar[ur]_V^{\ }="t2"
       \ar@{=>}|<<<{\mathrm{Id}} "s1"; "t2"
       \ar@{=>}|<<<{\mathrm{Id}} "s2"; "t1"
    }
    }
    \hspace{5pt}
    =
    \hspace{5pt}
    \xymatrix{
      \bullet
      \ar@/^1.7pc/[rr]^{V \otimes V}_{\ }="s"
      \ar@/_1.7pc/[rr]_{V \otimes V}^{\ }="t"
      &&
      \bullet
      \ar@{=>}^{\theta_{V,V}} "s"; "t"
    }
   \,,
$$ respectively, where $\xymatrix{V \otimes W \ar[r]^{\theta_{V,W}} & W \otimes
V}$ denotes the canonical symmetric braiding isomorphism in $\mathrm{Vect}$.

The monoids assigned by the corresponding local nets $\mathcal{A}_{Z_\Vert}$ and
$\mathcal{A}_{Z_\times}$ are algebras of the form $\mathrm{End}(V^{\otimes n})$,
where $n$ is the total number of elementary edges in the respective boundary of
a region.

Given the inclusion of regions $O_{a,b} \subset O_{x,x'}$ 
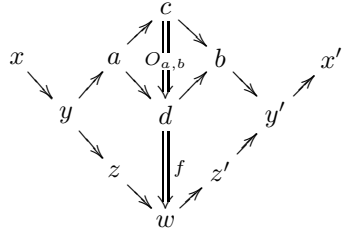
$$
  \xymatrix@R=7pt@C=7pt{
    &&&
    c
    \ar[dr]
    \ar@{=>}[dd]|{O_{a,b}}
    \\
    x
    \ar[dr]
    &&
    a
    \ar[ur]
    \ar[dr]
    &&
    b
    \ar[dr]
    &&
    x'
    &&&
    \\
    &
    y
    \ar[dr]
    \ar[ur]
    &&
     d
     \ar[ur]
     \ar@{=>}[dd]^f
    &&
    y'
    \ar[ur]
    &&&
    \\
    &&
    z
    \ar[dr]
    &&
    z'
    \ar[ur]
    &&&
    \\
    &&&
    w
    \ar[ur]
  }
$$ we get, according to definition \ref{nets from functors}, inclusions $$
  A_{Z_\Vert}, A_{Z_\times}
  :
  \mathrm{End}(V^{\otimes 2})
  \hookrightarrow
  \mathrm{End}(V^{\otimes 6})
$$ of endomorphism algebras given by $$
  A_{Z_\Vert}
  :
  \left(
    \begin{array}{cc}
      A & B
      \\
      C & D
    \end{array}
  \right)
  \mapsto
  \left(
    \begin{array}{cccccc}
      1 & 0 & 0 & 0 & 0 & 0
      \\
      0 & 1 & 0 & 0 & 0 & 0
      \\
      0 & 0 & A & B & 0 & 0
      \\
      0 & 0 & C & D & 0 & 0
      \\
      0 & 0 & 0 & 0 & 1 & 0
      \\
      0 & 0 & 0 & 0 & 0 & 1
    \end{array}
  \right)
  \,;
  \hspace{7pt}
  A_{Z_\times}
  :
  \left(
    \begin{array}{cc}
      A & B
      \\
      C & D
    \end{array}
  \right)
  \mapsto
  \left(
    \begin{array}{cccccc}
      1 & 0 & 0 & 0 & 0 & 0
      \\
      0 & A & 0 & 0 & B & 0
      \\
      0 & 0 & 1 & 0 & 0 & 0
      \\
      0 & 0 & 0 & 1 & 0 & 0
      \\
      0 & C & 0 & 0 & D & 0
      \\
      0 & 0 & 0 & 0 & 0 & 1
    \end{array}
  \right)
  \,,
$$ where each entry in these matrices is an endomorphism of $V$.

The locality of the net $\mathcal{A}_{Z_\Vert}$ is manifest. The algebras
assigned to two elementary regions clearly commute if and only if the two
regions are spacelike separated. For $\mathcal{A}_{Z_\times}$ the algebras of
course also commute if the regions are spacelike separated, but here they also
commute if the two regions are \emph{timelike} separated. Only if two elementary
regions are lightlike separated do the inclusions of algebras due to
$\mathcal{A}_{Z_\times}$ not commute.

There are various variations of this example. In particular for $Z_\times$ one
would want to consider the case where two different vector spaces $V_l$ and
$V_r$ and two nontrivial automorphisms $U_l : V_l \to V_l$ and $U_r : V_R \to
V_r$ are assigned to elementary causal subsets as follows: $$
  Z_\times
  \left(
    \raisebox{40pt}{
    \xymatrix{
       & y
       \ar[dr]
       \ar@{=>}[dd]
       \\
       x
       \ar[ur]
       \ar[dr]
       &&
       z
       \\
       &
       y'
       \ar[ur]
    }
    }
  \right)
  \hspace{7pt}
    :=
  \hspace{7pt}
    \raisebox{40pt}{
    \xymatrix{
       & \bullet
       \ar[dr]^{V_r}_{\ }="s2"
       \\
       \bullet
       \ar[ur]^{V_l}_{\ }="s1"
       \ar[dr]_{V_r}^{\ }="t1"
       &&
       \bullet
       \\
       &
       \bullet
       \ar[ur]_{V_l}^{\ }="t2"
       \ar@{=>}|<<<{U_l} "s1"; "t2"
       \ar@{=>}|<<<{U_r} "s2"; "t1"
    }
    }
    \hspace{5pt}
    =
    \hspace{5pt}
    \xymatrix{
      \bullet
      \ar@/^1.7pc/[rr]^{V_l \otimes V_r}_{\ }="s"
      \ar@/_1.7pc/[rr]_{V_r \otimes V_l}^{\ }="t"
      &&
      \bullet
      \ar@{=>}|{\theta_{V_l,V_r} \circ U_l \otimes U_r} "s"; "t"
    }
   \,,
$$ Denote by $$
  c : \mathrm{End}(V_r)\otimes \mathrm{End}(V_l)
   \hookrightarrow
   \mathrm{End}(V_r \otimes V_l)
$$ the canonical inclusion of algebras and by $$
  \xymatrix{
    c^* \mathcal{A}_{Z_\times}
    \ar@{^{(}->}
    [r]
    &
    \mathcal{A}_{Z_\times}
  }
$$ the local sub-net of $\mathcal{A}_{Z_\times}$ obtained by restricting along
$c$ everywhere. Then $c^*\mathcal{A}_{Z_\times}$ is what is called a
\emph{chiral} AQFT. Its structure is encoded entirely in the two independent
projections onto two orthogonal lightlike curves. $$
  c^* \mathcal{A}_{Z_\times}
  :
    \raisebox{40pt}{
    \xymatrix{
       & y
       \ar[dr]
       \ar@{=>}[dd]
       \\
       x
       \ar[ur]
       \ar[dr]
       &&
       z
       \\
       &
       y'
       \ar[ur]
    }
    }
    \hspace{3pt}
    \mapsto
    \hspace{3pt}
    \mathcal{A}_l
    \left(
    \raisebox{20pt}{
    \xymatrix{
       &
       z
       \\
       y'
       \ar[ur]
    }
    }
    \right)
    \otimes
    \mathcal{A}_r
    \left(
    \raisebox{20pt}{
    \xymatrix{
       x
       \ar[dr]
       \\
       &
       y'
    }
    }
    \right)
    =
    \mathrm{End}(V_l) \otimes \mathrm{End}(V_r)
    \,.
$$ Restricting attention to just one of these and then ``compactifying'' that to
a circle leads to the models \cite{KawahigashiI,KawahigashiLongo} of
2-dimensional (conformal) field theories as local nets on the circle.

This important example is further expanded on in section \ref{boundary}.

\subsection{Boundary FQFT and boundary AQFT} \label{boundary}

AQFT on spaces with boundary has been introduced in \cite{LongoRehren} for the
case of the Minkowski half-plane $X = \mathbb{R}^2_{<}$. Here we briefly
indicate how boundary conditions are formulated for FQFT and how we recover the
picture in \cite{LongoRehren} from this point of view.

We obtain the poset of causal subsets on the half plane, $S(\mathbb{R}^2_<)$, by
starting with $S(\mathbb{R}^2)$ and intersecting everything with
$\mathbb{R}^2_<$. We form $P_2(\mathbb{R}^2_<)$ by first restricting to 2-paths
that run entirely within $\mathbb{R}^2_<$ and then throwing in new boundary
generators for 1- and 2-morphisms of the form $$
    
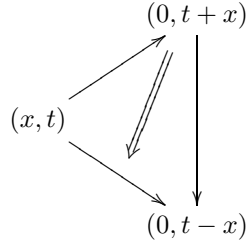
\xymatrix{
       & (0,t+x)
       \ar[dd]_{\ }="s2"
       \\
       (x,t)
       \ar[ur]_>{\ }="s1"
       \ar[dr]^{\ }="t1"
       \\
       &
       (0,t-x)
       \ar@{=>} "s1"; "t1"
    }
$$

From examples of classical parallel $n$-transport \cite{ndclecture} and from the
2-functorial description of rational CFT \cite{FS} it is known that boundary
conditions for $n$-functors $Z$ correspond to choices of morphism from some
trivial $n$-functor $I$ into the restriction of the given one to the boundary:
$$
  \xymatrix{
    I
    \ar[rr]
    &&
    Z|_{\partial X}
  }
  \,.
$$ We illustrate this in the context of the last example, $
  Z_\times : P_2(\mathbb{R}^2) \to \mathbf{B}\mathrm{Vect}
$, from section \ref{examples}, which lead to the discussion of chiral nets $i^*
\mathcal{A}_{Z_\times} \subset \mathcal{A}_{Z_\times}$.

For that purpose, let $I$ be the 2-functor $
  I : P_2(\mathbb{R}^2) \to \mathbf{B}\mathrm{Vect}
$ which is constant on the single object of $\mathbf{B}\mathrm{Vect}$ and
consider 2-functors $
  Z_\times^{<} :
  \xymatrix{
    P_2(\mathbb{R}^2_<) \ar[r]
    &
    \mathbf{B}\mathrm{Vect}
  }
$ which coincide with our $Z_\times$ in the bulk. Then we have the simple but
important \begin{proposition}
 If a morphism
 $$
   b : I \to Z^<_\times|_{\partial \mathbb{R}^2_<}
 $$
 exists and
 is time independent in that its component map is constant
 on objects (but not the 0 dimensional vector space), then
 $Z^<_\times$ assigns the identity to all boundary paths.
\end{proposition} \proof
  The components of the morphism, which is a pseudonatural transformation
  of 2-functors, are 2-cells in $\mathbf{B}\mathrm{Vect}$ of the form
  $$
    \raisebox{30pt}
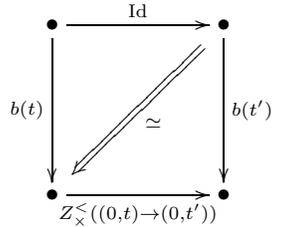
{
    \xymatrix{
      \bullet
      \ar[rr]^{\mathrm{Id}}
      \ar[dd]_{b(t)}^>{\ }="t"
      && \bullet
      \ar[dd]^{b(t')}_<{\ }="s"
      \\
      \\
      \bullet
      \ar[rr]_{Z^<_\times((0,t)\to(0,t'))}
      &&
      \bullet
      \ar@{=>}^\simeq "s"; "t"
    }
    }
    \,.
  $$
  By assumption of time independence of the boundary condition we have
  $b(t) = b(t') = b(0)$.
  This means that $Z^<_\times((0,t)\to(0,t'))$ must be a vector space
  such that there exists an isomorphism of vector spaces
  $$
     b(0)\otimes Z^<_\times((0,t)\to(0,t'))
     \simeq
     b(0)
     \,.
  $$
\endofproof

So in this case the 2-functor $Z^<_\times$ will specify identifications of the
vector spaces $V_l$ and $V_r$ at the boundary $$
  Z^<_\times
  \hspace{4pt}
  :
  \hspace{4pt}
  \raisebox{43pt}{
    \xymatrix{
       & (0,t+x)
       \ar[dd]_{\ }="s2"
       \\
       (x,t)
       \ar[ur]_>{\ }="s1"
       \ar[dr]^{\ }="t1"
       \\
       &
       (0,t-x)
       \ar@{=>} "s1"; "t1"
    }
    }
    \hspace{7pt}
     \mapsto
    \hspace{7pt}
    \raisebox{40pt}{
    \xymatrix{
       & \bullet
       \ar[dd]^{\mathrm{Id}}_{\ }="s2"
       \\
       \bullet
       \ar[ur]^{V_l}_{\ }="s1"
       \ar[dr]_{V_r}^{\ }="t1"
       \\
       &
       \bullet
       \ar@/^.6pc/@{=>}_<<<{\simeq} "s1"; "t1"
    }
   }
   \,.
$$

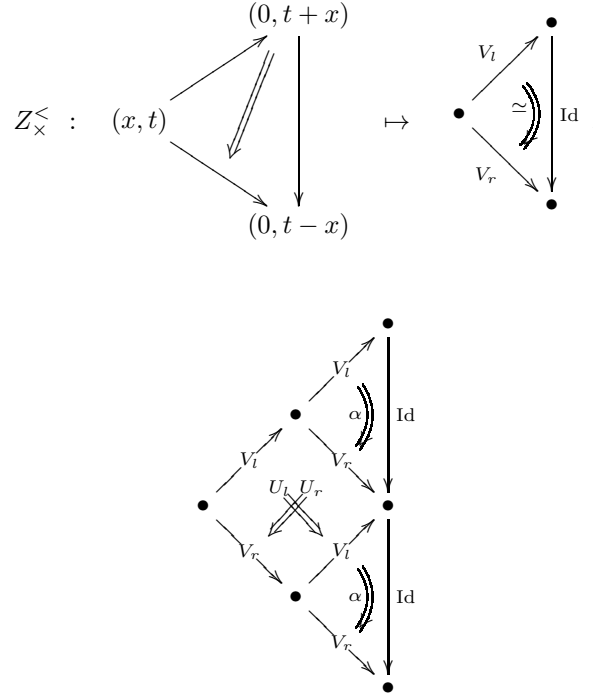
\begin{figure}[h] $$
  \xymatrix{
     &
     &\bullet
     \ar[dd]^{\mathrm{Id}}
     \\
     &
     \bullet
      \ar[ur]|{V_l}_{\ }="s3"
      \ar[dr]|{V_r}_{\ }="s2"^{\ }="t3"
     \\
     \bullet
     \ar[ur]|{V_l}_{\ }="s1"
     \ar[dr]|{V_r}^{\ }="t2"
     &
     & \bullet
     \ar[dd]^{\mathrm{Id}}
     \\
     &
     \bullet
     \ar[ur]|{V_l}^{\ }="t1"_{\ }="s4"
     \ar[dr]|{V_r}^{\ }="t4"
     \\
     &
     & \bullet
     \ar@{=>}|<<<{U_l} "s1"; "t1"
     \ar@{=>}|<<<{U_r} "s2"; "t2"
     \ar@/^.5pc/@{=>}_\alpha "s3"; "t3"
     \ar@/^.5pc/@{=>}_\alpha "s4"; "t4"
  }
$$ \caption{
  The image under the boundary FQFT 2-functor $Z^<_\times$ of a
  spacelike wedge on the left Minkowski half plane.
} \end{figure}

By taking endomorphisms this defines a net of algebras on the boundary, which
entirely encodes the chiral part $c^*\mathcal{A}_{Z^<_\times}$ of
$\mathcal{A}_{Z^<_\times}$. This way we arrive at the picture of boundary AQFT
given in \cite{LongoRehren}. Further details should be discussed elsewhere.

\subsection{2-$C^*$-category codomains}

In most applications to physics one wants the algebras in a local
net to be $C^*$-algebras. A natural type of 2-category in which
endomorphism algebras of 1-morphisms are $C^*$-algebras is that of
2-$C^*$-categories: categories enriched in $C^*$-categories.

\begin{definition}
  A \emph{$C^*$-category} (or \emph{$C^*$-algebroid}: the many-object
version of a $C^*$-algebra) is a category $C$ enriched in
complex Banach spaces (meaning that for all objects $\rho,\sigma,\tau$ of $C$
we have that $C(\rho,\sigma)$ is a complex Banach space and that composition
$$
  \circ_{\rho,\sigma,\tau} : C(\rho,\sigma)
  \times C(\sigma,\tau) \to C(\rho,\tau)
$$
is a morphism of complex Banach spaces)
which is equipped with an involutive antilinear functor
$$
  (\cdot)^* : C \to C^{\mathrm{op}}
$$
that satisfies the
    \emph{$C^*$-condition}
$$
  \forall \rho,\sigma \in \mathrm{Obj}(C):\,
  \forall S \in C(a,b)\,:\,
  \left\{
  \begin{array}{l}
    \mbox{$S^* \circ S$ is positive in $C(\rho,\rho)$}
    \\
    \Vert
      S^* \circ S
    \Vert
     =
    \Vert S\Vert^2
  \end{array}
  \right.\,,
$$
  where $\Vert \cdot \Vert : C(\rho,\sigma) \to \mathbb{C}$ is the
  Banach norm.
\end{definition}
A $C^*$-algebra $A$ is precisely the endomorphism algebra
of an object $\rho$ in a $C^*$-category, $A = C(\rho,\rho)$.
We write $\mathbf{B}A$ for the one object $C^*$-category
whose single endomorphism algebra is $A$.

$C^*$-categories form a strict monoidal 2-category
$(C^*\mathrm{Cat},\times)$ whose morphisms are
Banach space functors (continuous on each Hom-space).
Therefore one can enrich in $C^*$-categories themselves:

\begin{definition}
  A (strict) \emph{2-$C^*$-category} is a category
  enriched in $C^*\mathrm{Cat}$.
\end{definition}
A discussion of aspects of 2-$C^*$-categories can be found in \cite{Zito}.

The canonical example of a strict 2-$C^*$-category is
$\mathrm{Ampli}_{C^*} \subset \mathrm{Bimod}_{C^*}$, the 2-category whose objects are
unital $C^*$-algebras, whose morphisms are amplimorphisms between
these and whose 2-morphisms are intertwiners between those.
$\mathrm{Bimod}_{C^*}$ is very similar, but is not strict.
See \cite{LR} and section 2 of \cite{Zito}.

So we have
\begin{observation}
  For $Z : P_2(X) \to \mathcal{C}$ a transport 2-functor
  with values in a 2-$C^*$-category $\mathcal{C}$, the corresponding local net $A_Z$ is a net of $C^*$-algebras.
\end{observation}

\subsection{Hopf spin chain models}
\label{Hopf spin chain models}

Recall the description of lattice models with boundary from section
\ref{boundary}. Consider the extreme case where there is a left and
right boundary which are separated only by a single lattice spacing:
$$
  \raisebox{80pt}
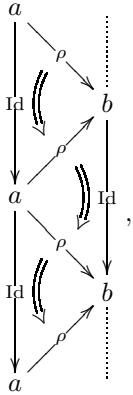
{
  \xymatrix{
    a
    \ar[dr]|\rho_{\ }="s1"
    \ar[dd]|{\mathrm{Id}}
    &
    \ar@{..}[d]
    \\
    &
    b
    \ar[dd]|{\mathrm{Id}}
    \\
    a
    \ar[ur]|{\rho}^{\ }="t1"_{\ }="s2"
    \ar[dd]|{\mathrm{Id}}
    \ar[dr]|{\rho}^{\ }="t2"_{\ }="s3"
    &
    \\
    & b
    \ar@{..}[d]
    \\
    a
    \ar[ur]|{\rho}^{\ }="t3"
    &
    \\
    \ar@/_.4pc/@{=>} "s1"; "t1"
    \ar@/^.4pc/@{=>} "s2"; "t2"
    \ar@/_.4pc/@{=>} "s3"; "t3"
  }
  }
  \,,
$$
where for simplicity we are concentrating on the case that
$Z$ sends each edge to one and the same morphism
$\rho : a \to b$ in $\mathcal{C}$.

Physically, we can think of this as a lattice model for an open string stretching from a brane of type $a$ to a brane of type $b$.
It s a crude lattice model, consisting of a single ``string bit''.

Consider another such strip, labeled by another morphism
$\bar \rho : b \to a$
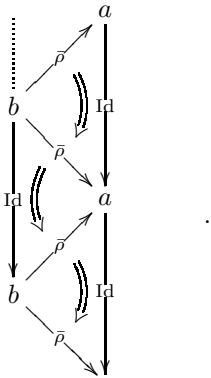
$$
  \raisebox{80pt}{
  \xymatrix{
    \ar@{..}[d]
    &
    a
    \ar[dd]|{\mathrm{Id}}
    &
    \\
    b
    \ar[ur]|{\bar \rho}_{\ }="s1"
    \ar[dr]|{\bar \rho}^{\ }="t1"_{\ }="s2"
    \ar[dd]|{\mathrm{Id}}
    &
    \\
    &
    a
    \ar[dd]|{\mathrm{Id}}
    \\
    b
    \ar[ur]|{\bar \rho}^{\ }="t2"_{\ }="s3"
    \ar[dr]|{\bar \rho}^{\ }="t3"
    &
    \\
    &
    \ar@/^.4pc/@{=>} "s1"; "t1"
    \ar@/_.4pc/@{=>} "s2"; "t2"
    \ar@/^.4pc/@{=>} "s3"; "t3"
  }
  }
  \,.
$$
As the notation suggests, we want to think of $\bar \rho$
to be \emph{conjugate} to $\rho$, meaning that
$\rho$ and $\bar \rho$
form an ambidextrous adjunction \cite{Lauda}
between $a$ and $b$ such that the unit of the left-handed adjunction is the ${}^*$-adjoint of the counit of the right-handed adjunction, and vice versa.
(see p. 8 of \cite{Zito}).

Then it makes sense to think of this as a lattice model for an open string, or rather a ``string bit'', as before, but now with that string taken to stretch from the $b$-type brane to the $a$-type brane.
We can then consider lattice models built from the above building blocks by gluing the above strip-wise 2-functors horizontally:
$$
  \raisebox{80pt}
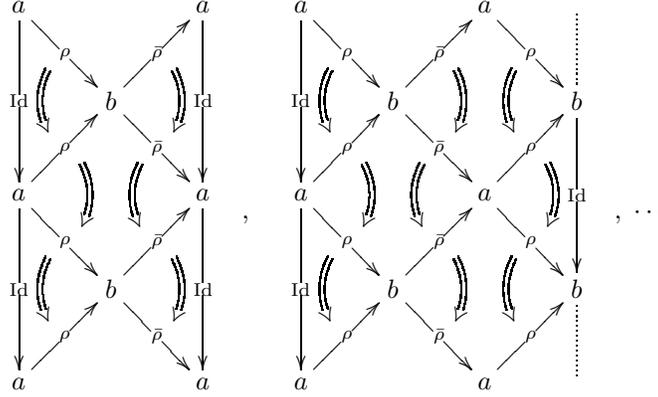
{
  \xymatrix{
    a
    \ar[dr]|\rho_{\ }="s1"
    \ar[dd]|{\mathrm{Id}}
    &
    &
    a
    \ar[dd]|{\mathrm{Id}}
    \\
    &
    b
    %\ar[dd]|{\mathrm{Id}}
    \ar[ur]|{\bar \rho}_{\ }="ss1"
    \ar[dr]|{\bar \rho}^{\ }="tt1"_{\ }="ss2"
    &
    \\
    a
    \ar[ur]|{\rho}^{\ }="t1"_{\ }="s2"
    \ar[dd]|{\mathrm{Id}}
    \ar[dr]|{\rho}^{\ }="t2"_{\ }="s3"
    &
    &
    a
    \ar[dd]|{\mathrm{Id}}
    \\
    & b
    %\ar@{..}[d]
    \ar[ur]|{\bar \rho}^{\ }="tt2"_{\ }="ss3"
    \ar[dr]|{\bar \rho}^{\ }="tt3"
    \\
    a
    \ar[ur]|{\rho}^{\ }="t3"
    &
    &
    a
    \\
    \ar@/_.4pc/@{=>} "s1"; "t1"
    \ar@/^.4pc/@{=>} "s2"; "t2"
    \ar@/_.4pc/@{=>} "s3"; "t3"
    \ar@/^.4pc/@{=>} "ss1"; "tt1"
    \ar@/_.4pc/@{=>} "ss2"; "tt2"
    \ar@/^.4pc/@{=>} "ss3"; "tt3"
  }
  }
  \hspace{9pt}
  ,
  \hspace{9pt}
  \raisebox{80pt}{
  \xymatrix{
    a
    \ar[dr]|\rho_{\ }="s1"
    \ar[dd]|{\mathrm{Id}}
    &
    &
    a
    %\ar[dd]|{\mathrm{Id}}
    \ar[dr]|{\rho}_{\ }="sss1"
    &
    \ar@{..}[d]
    \\
    &
    b
    %\ar[dd]|{\mathrm{Id}}
    \ar[ur]|{\bar \rho}_{\ }="ss1"
    \ar[dr]|{\bar \rho}^{\ }="tt1"_{\ }="ss2"
    &
    &
    b
    \ar[dd]|{\mathrm{Id}}
    \\
    a
    \ar[ur]|{\rho}^{\ }="t1"_{\ }="s2"
    \ar[dd]|{\mathrm{Id}}
    \ar[dr]|{\rho}^{\ }="t2"_{\ }="s3"
    &
    &
    a
    \ar[ur]|{\rho}^{\ }="ttt1"_{\ }="sss2"
    \ar[dr]|{\rho}^{\ }="ttt2"_{\ }="sss3"
    \\
    & b
    %\ar@{..}[d]
    \ar[ur]|{\bar \rho}^{\ }="tt2"_{\ }="ss3"
    \ar[dr]|{\bar \rho}^{\ }="tt3"
    &&
    b
    \ar@{..}[d]
    \\
    a
    \ar[ur]|{\rho}^{\ }="t3"
    &
    &
    a
    \ar[ur]|{\rho}^{\ }="ttt3"
    &
    \ar@/_.4pc/@{=>} "s1"; "t1"
    \ar@/^.4pc/@{=>} "s2"; "t2"
    \ar@/_.4pc/@{=>} "s3"; "t3"
    \ar@/^.4pc/@{=>} "ss1"; "tt1"
    \ar@/_.4pc/@{=>} "ss2"; "tt2"
    \ar@/^.4pc/@{=>} "ss3"; "tt3"
    \ar@/_.4pc/@{=>} "sss1"; "ttt1"
    \ar@/^.4pc/@{=>} "sss2"; "ttt2"
    \ar@/_.4pc/@{=>} "sss3"; "ttt3"
  }
  }
  \hspace{9pt}
  ,
  \hspace{3pt}
  \cdots
$$
The algebras assigned by the corresponding net $A_Z$ to the elementary causal bigon $O_{\rho,\bar \rho}$ and
$O_{\bar \rho,\rho}$ are
$$
  A_Z(O_{\rho,\bar \rho}) = \mathrm{End}_{\mathcal{C}}(\bar \rho \circ \rho)
$$
and
$$
  A_Z(O_{\bar \rho,\rho}) = \mathrm{End}_{\mathcal{C}}(\rho \circ \bar \rho)
  \,.
$$
If $\mathcal{C}$ is a 2-$C^*$-category, these are $C^*$-Hopf algebras $H$ and $\hat H$ which are duals of each other \cite{Mueger, Zito}.
Due to the fact that the 2-morphisms in the above diagrams do not mix $\rho$ and $\bar \rho$, we can understand the nature of the net $A_Z$ obtained from the above 2-functor $Z$ already by concentrating on the endomorphism algebras assigned to a horizontal zig-zag
$$
 \raisebox{20pt}
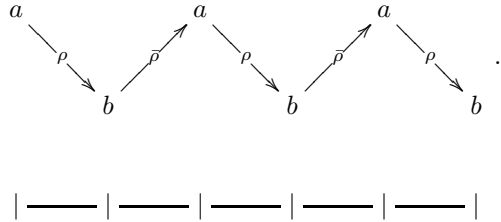
{
  \xymatrix{
    a
    \ar[dr]|{\rho}
    &&
    a
    \ar[dr]|{\rho}
    &&
    a
    \ar[dr]|{\rho}
    \\
    &
    b \ar[ur]|{\bar \rho}
    &
    &
    b \ar[ur]|{\bar \rho}
    &
    &
    b %\ar[ur]|{\bar \rho}
    \\
    |
    \ar@{-}[r]
    &
    |
    \ar@{-}[r]
    &
    |
    \ar@{-}[r]
    &
    |
    \ar@{-}[r]
    &
    |
    \ar@{-}[r]
    &
    |
  }
  }
  \,.
$$
If we to restrict to evaluating the net $A_Z$ on zig-zags of even length, this gives rise to a net on the latticized real axis
with the property that algebras $A_Z(I_1)$ and $A_Z(I_2)$
commute if the intervals $I_1$ and $I_2$ are not just disjoint but
differ by at least one lattice spacing.
Precisely these kind of 1-dimensional nets are considered in
\cite{NS}, where they are addressed as
\emph{Hopf spin chain models}.

\subsection{Subfactors and asymptotic inclusion}
\label{subfactors and asymptotic inclusion}

We can interpret the analysis of the direct limit algebras of the
above nets $A_Z$ given in \cite{NS} in terms of Ocneanu's notion of
asymptotic inclusion \cite{Ocneanu} and its relation to subfactors.
(Thanks to Pasqual Zito for pointing this out.)

Consider first a lattice of the above sort unbounded (only) to the
right. The direct limit algebra
$$
  A :=  \mathrm{colim}\,  A_Z
$$
of the chain of inclusions of finite algebras
$$
  \mathrm{End}_{\mathcal{C}}(\bar \rho \circ \rho)
  \hookrightarrow
  \mathrm{End}_{\mathcal{C}}(\rho \circ \bar \rho \circ \rho)
  \hookrightarrow
  \mathrm{End}_{\mathcal{C}}(\bar \rho \circ \rho \circ \bar \rho \circ \rho)
  \hookrightarrow
  \cdots
$$
naturally carries a trace, which we can assume to be normalized. Completing with respect to the norm $\Vert a \Vert := \mathrm{tr}(a^* \circ a)$ yields an algebra $\bar A$ which is a type II vonNeumann algebra factor.

We can shift everything one lattice spacing to the right and consider the poset of algebras
$$
  \mathrm{Id}_\rho \cdot \mathrm{End}_{\mathcal{C}}(\bar \rho)
  \hookrightarrow
  \mathrm{Id}_\rho \cdot
   \mathrm{End}_{\mathcal{C}}(\rho \circ \bar \rho )
  \hookrightarrow
  \mathrm{Id}_\rho \cdot  \mathrm{End}_{\mathcal{C}}(\bar \rho \circ \rho \circ \bar \rho )
  \hookrightarrow
  \cdots
  \,,
$$
where $\cdot$ denotes the horizontal composition in our
2-$C^*$-category $\mathcal{C}$. The completion of the
direct limit of this chain of inclusions is a type II factor
$\bar B$ which has a canonical inclusion into $\bar A$
$$
  \bar B \hookrightarrow \bar A
  \,.
$$
This inclusion of subfactors obtained from a pair of conjugate
morphisms $\rho, \bar \rho$ in a 2-$C^*$-category is Ocneanu's
\emph{asymptotic inclusion} \cite{Ocneanu,EK}. 
If the
2-$C^*$-category $\mathcal{C}$ that we started with is $\mathcal{C}
= \mathrm{Bimod}_{C^*}$ and the original morphism $\rho : a \to b$
in $\mathcal{C}$ itself an inclusion of subfactors, then this is
recovered by the above construction.

$\bar A$ and $\bar B$ and their inclusion
$\bar B \hookrightarrow \bar A$ encode a QFT on the right half plane. From the above setup we can analogously obtain a subfactor
$\bar B^o \hookrightarrow \bar A^o$ for the left half plane.
Moreover, the completion of the direct limit algebra over all endomorphism algebras of zig-zags that are allowed to extend finitely to the right and the left yields a factor $\bar K$ which has a canonical inclusion of the factor $\bar A^o \otimes \bar A$
$$
  \bar A^o \otimes \bar A
  \hookrightarrow
  K
  \,.
$$
Following the
discussion on p. 10 of \cite{Kawahigashi} one can understand this in
the context of \cite{LR} and read $\bar A$ and $\bar A^o$ as two
chiral open string algebras and $K$ as the corresponding closed
string algebra.

\section{Further issues} \label{further issues}

There are various immediate further questions to be addressed. We shall be
content here with just briefly commenting on the following four.

\subsection{General Lorentzian structure} AQFT was originally conceived
entirely in its application to quantum field theories on Minkowski space, which
is the case we have been concentrating on above. A generalization of
Poincar{\'e}-covariant nets on causal subsets in Minkowski space to nets on
globally hyperbolic Lorentian spaces has later been proposed in
\cite{generallycovariant}.

The possibly most natural and immediate generalization to AQFT on a fixed
general Lorentzian space was indicated in \cite{locallocality}: on a
Lorentzian manifold $X$ an AQFT net should be \emph{locally local}: the
locality axiom should hold after restriction of the net to any globally
hyperbolic subspace of $X$. The same should be true for the time slice axiom.

No guesswork is required for generalizing the concept of Minkowskian FQFT
2-functors to general Lorentzian 2-functors: the concept of the 2-functor
itself makes unambiguous sense for any choice of 2-path 2-category in $X$. So we
can use our construction of local nets from 2-functors to \emph{derive} locality
properties of nets on Lorentzian spaces. Doing so confirms the idea of
\cite{locallocality}:

  Let $(X,g)$ be any 2-dimensional oriented and time-oriented
  Lorentzian manifold.

In generalization of definition \ref{causal subsets} consider
\begin{definition}
  A causal subset $O \subset X$ is a subset of a globally hyperbolic
  subset of $X$ which is the interior of a non-empty
  intersection of the future of
  one point with the past of another. Write $S(X)$ for the
  category with such causal subsets as objects and inclusion
  of subsets as morphisms.
\end{definition}

In generalization of definition \ref{2paths in Minkowski} consider
\begin{definition}
  Let $P_2(X)$ be the strict 2-category whose objects are the points in $X$,
  whose 1-morphisms are piecewise lightlike and right-moving paths
  (with respect to the chosen orientation and time-orientation of $X$)
  and whose 2-morphisms are
  generated under gluing along common boundaries from closures of
  causal subsets.
\end{definition}

Our construction in definition \ref{nets from functors} immediately generalizes
to a construction of a net $\mathcal{A}_Z : S(X) \to \mathrm{Monoids}$ from a
2-functor $Z : P_2(X) \to C$. All the arguments need to be done within globally
hyperbolic subsets of $X$, where they go through literally as before. We can
\emph{read off} from the result of this construction the locality properties of
$\mathcal{A}_Z$:

\begin{proposition}
  The net $\mathcal{A}_Z : S(X) \to \mathrm{Monoids}$
 obtained from any 2-functor $Z : P_2(X) \to C$ is \emph{locally local}
 and satisfies the \emph{local time slice axiom}:
 for any inclusion
 $$
   i : Y \hookrightarrow X
 $$
 with $Y$ globally hyperbolic we have that
 $i^* \mathcal{A}_Z$ is a local net satisfying the time slice axiom.

\end{proposition} This concept of local locality is compatible with
\cite{generallycovariant} but does not presuppose any covariance condition on
the net.

\subsection{Higher dimensional QFT}
\label{higher dimensional QFT}

We had considered, for ease of discussion, in definition \ref{exFQFT} the
2-category $P_2(X)$ whose 2-morphisms are generated from gluing the closures of
2-dimensional causal subsets along common boundaries. But nothing in our
constructions crucially depends on gluing of causal subsets, and in fact gluing
of causal subsets becomes less naural in higher dimensions. As the examples we
presented in section \ref{examples}, where we obtained FQFT 2-functors by
\emph{restricting} 2-functors on a larger 2-category of 2-paths to $P_2(X)$,
clearly indicate, the 2-category $P_2(X)$ can be replaced by any 2-category of
2-paths in $X$ which is large enough that every causal subset in $X$ can be
regarded as a 2-morphisms in there, so that every FQFT 2-functor can be
evaluated on causal subsets. And this statement then immediately generalizes to
higher dimensions.

For $X$ a $d$-dimensional Lorentzian manifold, we should take the
category $S(X)$ to be that whose objects are causal subsets in $X$, which are
are those subsets that arise within any globally hyperbolic subset of $X$ as the
interior of the future of one point with the past of another point. Morphisms
are inclusions.

The $d$-category $P_d(X)$ used to described Lorentzian FQFT on $X$ can be
any sub-$d$-groupoid of the path $d$-groupoid \cite{ndclecture} which is large
enough so that every causal subset in $X$ comes from a $d$-morphism in $P_d(X)$
and such that the obvious higher dimensional generalizations of the diagrams in
section \ref{aqft from fqft} exist in $P_d(X)$. In particular, one can always
use the \emph{full} path $d$-groupoid.

\begin{figure}[h]
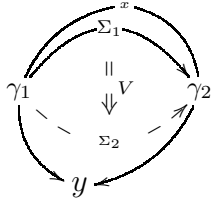
 $$ \xy
  (-12,0)*{\makebox(11,11){$\gamma_1$}}="left";
  (12,0)*{\makebox(10,11){$\gamma_2$}}="right";
  (2,11.5)*{\makebox(8,8){\tiny $x$}}="back";
  (-4,-12.8)*{\makebox(12,12){\large $y$}}="front";
  \ar@/^2pc/|{\Sigma_1} "left"; "right"
  \ar@/_1.6pc/@{-->}|{\mbox{\tiny $\Sigma_2$}} "left"; "right"
  \ar@/_.65pc/@{-} "back"; "left"
  \ar@/^.7pc/@{-} "back"; "right"
  \ar@/_.7pc/ "left"; "front"
  \ar@/^.7pc/ "right"; "front"
  \ar@{==>}^V (0,4); (0,-3)
\endxy $$
  \caption{
    A 3-morphism in a 3-path 3-category: a volume $V$, cobounding
    two surfaces $\Sigma_1$ and $\Sigma_2$, which each cobound
    two paths $\gamma_1$ and $\gamma_2$ which each cobound two
    points $x$ any $y$.
  }
\end{figure}

With such a setup, all our constructions here should have essentially
straightforward generalizations to higher dimensions, leading to a construction
of local nets on $X$ from any FQFT $d$-functor on $X$.

\subsection{Extended FQFT from AQFT?} We have shown how to go from FQFTs to
AQFTs. An obvious question is if there is a way to go back from AQFTs to FQFTs.
One would have to identify from a local net first of all the objects that the
local algebras are the endomorphism algebras of. Since these algebras are
usually $C^*$-algebras, this would be accomplished by using the Gelfand-Naimark
theorem, which states that every $C^*$-algebra is isomorphic to the
$C^*$-algebra of bounded operators on some Hilbert space. But to get a full
2-functorial FQFT, one needs also a compatible horizontal composition on these
Hilbert spaces. Potentially this can be extracted using the machinery of
\emph{localized transportable endomorphisms} as in section 8 of
\cite{HalvorsonMueger}.

\appendix

\section{2-Vector spaces and the canonical 2-representation} \label{2veccanrep}

In section \ref{examples} we obtained examples of FQFT 2-functors from
differential form data and a choice of 2-representation. Here we briefly
indicate a bit of background concerning these 2-representations.

For our purposes here a 2-vector space is an abelian module category, i.e an
abelian category equipped with an action by a monoidal category. Notice that the
category of $k$-vector spaces is the category of $k$-modules $$
  \mathrm{Vect}_k = k-\mathrm{Mod}
  \,.
$$ Accordingly we write $$
  2\mathrm{Vect} = \mathrm{Vect}_{\mathrm{Vect}} = \mathrm{Vect}-\mathrm{Mod}
$$ for the 2-category of abelian categories equipped with a (left, say)
$(\mathrm{Vect},\otimes)$-action. Since $\mathrm{Vect}$ is symmetric monoidal,
one can keep going this way and in principle define recursively the $n$-category
$$
  n\mathrm{Vect} = (n-1)\mathrm{Vect}-\mathrm{Mod}
  \,.
$$ Notice in particular that then $0\mathrm{Vect} = k$.

There are other monoidal categories over which one may want to consider 2-vector
spaces. For instance if we denote by $\mathrm{Disc}(k)$ the discrete category
over the ground field (the ground field as its objects and only identity
morphisms), then $$
  \mathrm{Disc}(k)-\mathrm{Mod} \simeq
  \mathrm{Cat}(\mathrm{Vect})
$$ is the 2-category of categories internal to vector spaces, which in turn is
equivalent to chain complexes concentrated in degree 0 and 1. These are the
2-vector spaces considered in \cite{BC}. $\mathrm{Disc}(k)$-modules are the
``right'' notion for 2-vector space for higher Lie theory, but probably not
\cite{BBK} as models for fibers of interesting 2-vector bundles.

The entirety of the 2-category of all $\mathrm{Vect}$-modules is quite
untractable. What is more accessible and more useful is the 2-category of
2-vector space that ``have a basis''. Noticing that an ordinary vector space $V$
has a basis if there is a set $S$ such that $V \simeq
\mathrm{Hom}_{\mathrm{Set}}(S,k)$, we should define a basis for a 2-vector space
$V$ to be a category $S$ such that $V \simeq \mathrm{Hom}(S,\mathrm{Vect})$. If
$S$ is itself $\mathrm{Vect}$-enriched this says that $V$ is a category of
algebroid modules. We shall restrict attention to $S$ having a single object, in
which case we are left with modules for ordinary algebras.

This way we find the bicategory $\mathrm{Bimod}$ of algebras, bimodules and
bimodule homomorphisms sitting inside $2\mathrm{Vect}$ as a sub-2-category of
2-vector spaces with basis: $$
  \xymatrix{
    \mathrm{Bimod}
    \ar@{^{(}->}[r]
    &
    2\mathrm{Vect}
  }
$$ $$
  \xymatrix{
    A
    \ar@/^2pc/[rr]^{N}_{\ }="s"
    \ar@/_2pc/[rr]_{N'}^{\ }="t"
    &&
    B
    \ar@{=>}^\phi "s"; "t"
  }
  \hspace{4pt}
    \mapsto
  \hspace{4pt}
  \xymatrix{
    \mathrm{Mod}_A
    \ar@/^2pc/[rr]^{-\otimes_A N}_{\ }="s"
    \ar@/_2pc/[rr]_{- \otimes_A N'}^{\ }="t"
    &&
    \mathrm{Mod}_B
    \ar@{=>}^{-\otimes_A \phi} "s"; "t"
  }
  \,.
$$ Notice how $\mathrm{Mod}_A$ is a \emph{category of modules} which is itself a
\emph{module category} over $\mathrm{Vect}$. The 2-category of
Kapranov-Voevodsky 2-vector spaces \cite{KV} is the full sub 2-category of
$\mathrm{Bimod}$ on all algebras of the form $k^{\oplus n}$ for $n \in
\mathbb{N}$. $$
 \mathrm{KV}2\mathrm{Vect}
  \hookrightarrow
  \mathrm{Bimod}
  \,.
$$

While $\mathrm{Bimod}$ is not a strict 2-category, it is a \emph{framed
bicategory} in the sense of \cite{Shulman}: there is the strict 2-category
$\mathrm{Algebras}$ of algebras, algebra homomorphisms and intertwiners (the
obvious 2-category for algebras regarded as one-object $\mathrm{Vect}$-enriched
categories), and the obvious inclusion $$
  \xymatrix{
    \mathrm{Algebras}
    \ar@{^{(}->}[r]
    &
    \mathrm{Bimod}
  }
$$ is full and faithful on all Hom-categories. Noticing that similarly groups,
when regarded as one-object groupoids, live in the 2-category $\mathrm{Groups}$
of groups, group homomorphisms and inner automorphisms, we get a strict
2-functor $$
  \xymatrix{
    \mathrm{Groups}
    \ar[r]
    &
    \mathrm{Algebras}
  }
$$ induced from forming for each group its group algebra. For each group $H$
there is the 2-group $\mathrm{AUT}(H) := \mathrm{Aut}_{\mathrm{Groups}}(H)$ and
the canonical inclusion $$
  \xymatrix{
  \mathbf{B}\mathrm{AUT}(H)
  \ar@{^{(}->}[r]
  &
  \mathrm{Groups}
  }
$$ induces, combined with the above discussion, the canonical 2-representation
of $\mathrm{AUT}(H)$ given by $$
  \rho_{\mathrm{can}}
  :
  \xymatrix{
    \mathbf{B}\mathrm{AUT}(H)
    \ar[r]
    &
    \mathrm{Groups}
    \ar[r]
    &
    \mathrm{Algebras}
    \ar[r]
    &
    \mathrm{Bimod}
    \ar[r]
    &
    2\mathrm{Vect}
  }
  \,.
$$ The logic of this construction generalizes to arbitrary strict 2-groups
$G_{(2)}$ coming from crossed modules of groups $(H \stackrel{t}{\to} G
\stackrel{\alpha}{\to} \mathrm{Aut}(G))$ (see for instance \cite{SWII} for a
review) and algebras obtained from a representation of $H$: \begin{proposition}
  For $\rho : \mathbf{B}H \to \mathrm{Vect}$ a representation of $H$
  such that the action of $G$ on $H$ extends to algebra automorphisms
  of the representation algebra $\langle \rho(H)\rangle$, the assignment
  $$
    \tilde \rho
    :
    \mathbf{B}(H \to G)
    \to
    \mathrm{Algebras}
  $$
  given by
  $$
  \xymatrix{
    \bullet
    \ar@/^2pc/[rr]^{g}_{\ }="s"
    \ar@/_2pc/[rr]_{g'}^{\ }="t"
    &&
    \bullet
    \ar@{=>}^h "s"; "t"
  }
  \mapsto
  \xymatrix{
    \langle \rho(H)\rangle
    \ar@/^2pc/[rr]^{\alpha(g)}_{\ }="s"
    \ar@/_2pc/[rr]_{\alpha(g')}^{\ }="t"
    &&
    \langle \rho(H)\rangle
    \ar@{=>}^{\rho(h)} "s"; "t"
  }
  $$
  is a strict 2-functor.
\end{proposition} Accordingly we obtain a 2-representation $$
 \xymatrix{
   \mathbf{B}(H \to G)
   \ar[r]^{\tilde \rho}
   &
   \mathrm{Algebras}
   \ar[r]
   &
   \mathrm{Bimod}
   \ar[r]
   &
   2\mathrm{Vect}
  }
   \,.
$$ All this should go through when the vector spaces here are equipped with more
structure. In particular, for $G$ a compact, simple and simply connected group,
for $\rho : \mathbf{B}\hat \Omega G \to \mathrm{Hilb}$ a positive-energy
representation of the weight 1 central extension of its loop group and for
$\mathrm{vNBimod}$ the bicategory of vonNeumann algebras and their bimodules
composed under Connes-fusion, \cite{ST} the above should extend to a
2-representation $$
  \mathbf{B}\mathrm{String}(G) \to
  \mathrm{vNBimod}
$$ of the strict String 2-group \cite{BCSS}.

\end{document}